\documentclass[11pt]{article}
\usepackage{graphics, latexsym,amssymb,amsmath,amscd}
\usepackage{graphics, latexsym,amssymb,amsmath,verbatim}
\usepackage{bm}
\date{January 3, 2022}
\title{Nonlinear Spencer operators on differentiable groupoids}
\author{ Jose M. M. Veloso\\  Faculdade de Matem\' atica - ICEN\\
Universidade Federal do Par\'a\\66059 - Bel\' em- PA - Brazil}
\begin{document}
\maketitle
\newcommand{\D}{\mbox{$\cal D$}}
\newtheorem{df}{Definition}[section]
\newtheorem{te}{Theorem}[section]
\newtheorem{co}{Corollary}[section]
\newtheorem{po}{Proposition}[section]
\newtheorem{lem}{Lemma}[section]
\newcommand{\zerobra}[2]{{{[}#1,#2{]}}}
\newcommand{\firstbra}[2]{{{[\![}#1,#2{]\!]}}}
\newcommand{\secondbra}[2]{{\rule[-3.03pt]{2pt}{12.05pt}{\hskip-2pt}{\bm{[}}{#1},{#2}{\bm{]}{\hskip-2pt}\rule[-3.03pt]{2pt}{12.05pt}}}}
\newcommand{\thirdbra}[2]{{\rule[-3.03pt]{3.5pt}{12.05pt}{\hskip-2pt}{\bm{[}}{#1},{#2}{\bm{]}{\hskip-2pt}\rule[-3.03pt]{3.5pt}{12.05pt}}}}

\newcommand{\Ad}{\mbox{Ad}}
\newcommand{\ad}{\mbox{ad}}
\newcommand{\im}[1]{\mbox{\rm im\,$#1$}}
\newcommand{\sime}{\mbox{sim}}
\newcommand{\Lg}{\mbox{$\frak g$}}
\newcommand{\Pf}{{\sc Proof}. }
\newcommand{\Ex}{{\sc Example}. }
\newcommand{\EPf}{\hbox{}\hfill$\Box$\vspace{.5cm}}
\newenvironment{proof}{{\normalsize {\sc Proof}:}}{{{\hfill $\Box$\vspace{.5cm}}}}

\begin{abstract}
We construct the first, second and sophisticated non-linear and linear Spencer complexes for a differentiable Lie groupoid $G$. To do this, we extend the diagonal calculus, as applied by Malgrange to the groupoid $M\times M$, to the context of $I\times G$, where $I$ is the manifold of identities of $G$.
\end{abstract}

\section{Introduction}
In \cite{S}, D. C. Spencer studying deformation of structures defined by transitive pseudo groups, introduced a nonlinear complex associated to the pseudogroup. There was work in the following years to simplify the introduction of this complex, as in Van Qu\^e \cite{N}, Kumpera and Spencer \cite{KS} and B. Malgrange \cite{Ma1}, \cite{Ma2}. The same ideas in the context of G-structures were developed by Guillemin and Sternberg \cite{GuS}. Malgrange obtained a construction for the nonlinear Spencer complex through the diagonal calculus of Grothendieck for the groupoid $M\times M$, where $M$ is a differentiable manifold. A open problem is to define the nonlinear Spencer complex for the groupoid $G_k$ of k-jets of bisections of a differentiable groupoid $G$ (see the Appendix of \cite{KS}). Our goal in this paper is to define this complex and introduce its properties. For this, we explore the relationship between vector fields on $I\times G$ and the actions of bisections of $G_{k+1}$ on sections of the algebroid  $J^k\mathfrak g$ associated to the groupoid $G_k$. Here $I$ is the submanifold of identities of $G$ and $\mathfrak g$ is the algebroid associated to $G$.

B. Malgrange in  \cite{Ma1}, \cite{Ma2}  considered the groupoid $M\times M$ and  $Q^k(M)$ the groupoid of k-jets of bisections of $M\times M$, and the algebroid $J^kTM$  associated to $Q^k(M)$. A section of $J^kTM$  is identified to  the quotient of a $\pi_1$ vertical vector field on $M\times M$ module the vector fields that are null up to order $k$ on the diagonal of $M\times M$. This identification is possible because $M\times M$ is a transitive groupoid, and a right invariant vector field on the %%source
$s$-fiber, $s$ the source of $G$, extends uniquely as a right invariant vector field all over the groupoid $M\times M$. When the groupoid is not transitive, the knowing of a right invariant vector field on the $s$-fiber is insufficient to extend it to all the groupoid $G$. To utilize the technique devised by Malgrange in the intransitive case we must consider right invariant vector fields on $I\times G$ by the right action of $G$ on $I\times G$ given by $(x,X).Y=(x,X.Y)$.  The diagonal of $M\times M$ is replaced by the diagonal $\Delta\subset I\times I\subset I\times G$. Let be $\mathcal R$ the sheaf of germs of vector fields on $I\times G$ that are $\rho_1$ projectables and right invariants ($\rho_1$ and $\rho_2$ are the projections of $I\times G$ on $I$ and $G$ respectively).  A section of $J^k\mathfrak g$ identifies with a $\rho_1$ vertical field of $\mathcal R$ module the vector fields in $\mathcal R$ that are null up to order $k$ on $\Delta$. The $\rho_2$ vertical vector fields in $\mathcal R$ identifies with $\mathcal T$ the sheaf of sections of the tangent space $T=TI$ of $I$. Thus $ \check{\mathcal J}^k\mathfrak g=\mathcal T\oplus {\mathcal J}^k\mathfrak g$ identifies to the sheaf $\mathcal R$ module the sub sheaf of $\mathcal R$ that are null up to order $k$ on $\Delta$. In \cite{Ma1}  $ \check{\mathcal J}^k\mathfrak g$ is obtained as sum of   ${\mathcal J}^k\mathfrak g$ and $\tilde{\mathcal J}^k\mathfrak g$. Here $\tilde{\mathcal J}^k\mathfrak g$ is the quotient of vector fields in $\mathcal R$  tangents to the submanifold  $\{(t(X),X)\in I\times G|X\in G\}$, where $t$ is the target of $G$. The action of bisections of $G_{k+1}$ on $\mathcal T$ was obtained through the actions on ${\mathcal J}^k\mathfrak g$ and $\tilde{\mathcal J}^k\mathfrak g$. When $G$ is intransitive, $\check{\mathcal J}^k\mathfrak g$ is not this sum. It is necessary to do a direct calculus to obtain the action of bisections of $G_{k+1}$ on $\mathcal T$ which is done in Proposition \ref{actionsigma}.

For another approach to Lie groupoids and algebroids, and Spencer operators, see \cite{CSS}. %Unfortunately in this presentation it is still missing a treatment for nonlinear Spencer operators and complexes. 
A more recent version of Malgrange construction is in \cite{V}.

We resume briefly the content of each section. In section \ref{p} we introduce the basic definitions of groupoids and algebroids, and the actions of jets of admissible sections on the jets of sections of tangent spaces. Section \ref{calcdiag} is the core section, where we introduce the diagonal calculus and the Lie algebra sheaf $\wedge( \check{\mathcal J}^\infty\mathfrak{g})^*\otimes ( \check{\cal J}^\infty\mathfrak g)$. In section \ref{fnsc} we introduce the first linear and non-linear Spencer complexes, through the introduction of sub-sheaf $\wedge\mathcal T^*\otimes \cal J^{\infty}\mathfrak g$, and give the basic properties of these complexes. In a similar way, in section \ref{snsc} we introduce the second linear and non-linear Spencer complexes, through the introduction of sub-sheaf $\wedge\tilde{\mathcal T}^*\otimes \tilde{\cal J}^{\infty}\mathfrak g$.  Finally, in section \ref{sSc}, we introduce the sophisticated linear and non linear Spencer complexes with their properties.

%I  thanks to M.A.Salazar by several discussions on this subject, which helped correct some imperfections in this work. Of course, errors are the responsibility of the author. 

This paper is dedicated to the memory of Alexandre Martins Rodrigues.

\section{Preliminaires}\label {p}

\subsection{Groupoids and algebroids}
\begin{df}
A \emph{differentiable groupoid} $G$ is a differentiable manifold $G$ with a regular submanifold $I$, two submersions $s, t:G\rightarrow I$ with $s^2=s$, $t^2=t$, and the following operations on $G$:
\begin{enumerate}
\item There exists a differentiable operation called \emph{composition} in $G$,
$$
\begin{array}{rcl}
(s\times t)^{-1}(\Delta)&\rightarrow &G,\\
(Y,X)&\rightarrow &YX
\end{array}
$$
where $\Delta=\{(x,x):x\in I\}$ with the following properties:
\begin{enumerate}
\item if $(Y,X),(Z,Y)\in (s\times t)^{-1}(\Delta)$ then $(Z,YX),(ZY,X)\in (s\times t)^{-1}(\Delta)$ and
$
Z(YX)=(ZY)X;
$
\item $Xs(X)=t(X)X=X.$
\end{enumerate}
\item There exists a diffeomorphism $\iota$ of $G$ called \emph{inversion} 
$$
\begin{array}{rcl}
\iota : G&\rightarrow &G,\\
X&\rightarrow &X^{-1}
\end{array}
$$
such that $(X,X^{-1}),(X^{-1},X)\in (s\times t)^{-1}(\Delta)$, and
$$
X^{-1}X=s(X), XX^{-1}=t(X).
$$

\end{enumerate}
The projections $s$ and $t$ are called of \emph{source} and \emph{target} respectively.
\end{df}

\Ex If $M$ is a differentiable manifold , $M\times M$ is a differentiable groupoid with $s=(\pi_2,\pi_2)$, $t=(\pi_1,\pi_1)$ and operations 
$(z,y)(y,x)=(z,x)$ and $(y,x)^{-1}=(x,y)$.

\Ex Let $(E,M,\pi)$ be a differentiable vector bundle  and $P$ the set of linear isomorphisms between the fibers  of $E$. This means that $X\in P$ if $X:E_a\rightarrow E_b$
is a linear isomorphism. The composition and inversion in $P$ are the composition and inversion of linear transformations. The identities are the identities $I_a$ in each $E_a$. Therefore the manifold $I$ of identities is diffeomorphic to $M$. The source map $s$ is defined by $s(X)=I_a$ and the target map $t$ by $t(X)=I_b$ for $X:E_a\rightarrow E_b$.

We denote by $G(x)=s^{-1}(x)$ the $s$-fiber of $G$ on $x\in I$; by $G(\cdot,y)=t^{-1}(y)$ the $t$-fiber of $G$ on $y\in M$ and $G(x,y)=G(x)\cap G(\cdot,y)$. The set $G(x,x)$ is a group, the so called \emph{isotropy group} of $G$ at point $x$. If $U, V$ are  open sets of $I$, we introduce the notations $G(U)=\cup_{x\in U}G(x)$, $G(.,V)=\cup_{y\in V}G(.,y)$, and $G(U,V)=G(U)\cap G(.,V)$.

A \emph{(differentiable) section} $F$ of $G$ defined on an open set $U$ of $M$ is a differentiable map $F:U\rightarrow G$ such that $s (F(x))=x$. If $t (F(U))=V$ and $f=t\circ F:U\rightarrow V$ is a diffeomorphism, we say that the section $F$ is a \emph{bisection}. We write $U=s(F)$, $V=t(F)$ and $t\circ F=tF$.

We denote by $\mathcal{G}$ the set of  bisections of $G$. Naturally $\mathcal{G}$ has a structure of groupoid.  If $F,H\in\mathcal{G}$ with $t (F)=s (H)$, then $HF(x)=H(f(x))F(x)$ and $F^{-1}(y)=F(f^{-1}(y))^{-1}$, $y\in t (F)$ where $f=tF$.

\subsection{Actions on $TG$}
A bisection $F$, with $s(F)=U$, $t(F)=V$, defines a diffeomorphism 
$$\begin{array}{rcl}
	\tilde F:G(\cdot,U)&\rightarrow & G(\cdot,V)\\
	X&\mapsto &F( t(X))X.
	\end{array}
	$$
The differential $\tilde F_*: TG(\cdot,U)\rightarrow TG(\cdot,V)$ depends, for each $X\in G(\cdot,V)$, only of $j^1_{t(X)}F$. This defines an action
\begin{equation}\label{jetf}
\begin{array}{rcl}
	j^1_{t(X)}F:T_XG&\rightarrow & T_{F(t(X))X}G\\
	v&\mapsto &j^1_{t(X)}F\cdot v=(\tilde F_*)_X(v).
	\end{array}
\end{equation}
The application (\ref{jetf}) defines a left action of the set $G_1$ of 1-jets of bisections of $G$ on $TG$
\begin{equation}\label{actionl}
\begin{array}{rcl}
	G_1\times TG&\rightarrow & TG\vspace{.2cm}\\
	(j^1_{t(X)}F,v\in T_XG)&\mapsto &j^1_{t(X)}F\cdot v\in T_{F(t(X))X}G.
	\end{array}
\end{equation} 
If $V_t\subset TG$ denotes the sub vector bundle of $t_*$ vertical vectors, then the action (\ref{actionl}) depends only on 
$F(t(X))$:
\begin{equation*}\label{restactionl}
\begin{array}{rcl}
	G\times V_{t}&\rightarrow & V_{t}\vspace{.2cm}\\
	(F({t(X)}),v\in (V_{t})_X)&\mapsto &F(t(X))\cdot v\in (V_{t})_{F(t(X))X}.
	\end{array}
\end{equation*} 
In a similar way, $F$ defines a right action which is a diffeomorphism
\begin{equation}\label{barF}
\begin{array}{rcl}
	\overline F:G(V)&\rightarrow & G(U)\\
	X&\mapsto &XF(f^{-1}(s(X))).
	\end{array}
\end{equation}
The differential $\overline F_*$ of $\overline F$ induces the right action
\begin{equation}\label{actionr}
\begin{array}{rcl}
	TG\times G_1&\rightarrow & TG\vspace{.2cm}\\
	(v\in T_XQ^k,j^1_{f^{-1}(s(X))}F)&\mapsto &v\cdot j^1_{f^{-1}(s(X))}F=(\overline F_*)_X(v).
	\end{array}
\end{equation}
As $t(YX)=t(Y)$, it follows that
\begin{equation*}\label{f15}
t_*(v\cdot j^1_{f^{-1}(s(X))}F)=t_*(v),
\end{equation*}
where $t_*:TG\rightarrow TI$ is the differential of $t:G\rightarrow I$.
We deduce from (\ref{barF}) that the function $\overline F$ restricted to the $s$-fiber $G(y)$ depends only on the value of $F$ in $f^{-1}(y)$. If $VG=\ker s_*$, then the right action (\ref{actionr}) depends only on the value of $F$ at each point, and %the action (\ref{actionr}) 
by restriction gives the action
\begin{equation}\label{actionvert}
\begin{array}{rcl}
	VG\times G&\rightarrow & VG\\
	(v\in V_XG,Y\in G(.,s(X))&\mapsto &v\cdot Y\in V_{XY}G.
	\end{array}
\end{equation}

A vector field $\overline\xi$
on $G$ with values in $VG$ is said \emph{right invariant} if $\overline\xi(XY)=\overline\xi(X)\cdot Y$. The vector field $\overline\xi$
is determined by its restriction $\xi$ to $I$. 

If $\xi$ is a section of $V$ on $U\subset I$, let be  
$$
\overline\xi(X)=\xi(t(X))\cdot X,
$$ 
the right invariant vector field on $G(.,U)$. Then $\overline\xi$ has $\overline F_u$, $-\epsilon<u<\epsilon$, as one parameter group of diffeomorphisms induced by  bisections $F_u$ of $G$ such that 
$$\frac{d}{du}\overline F_u|_{u=0}=\overline\xi.
$$
Therefore, $F_0=I$ and 
$$\frac{d}{du} F_u(x)|_{u=0}=\xi(x).
$$

\begin{df}
The vector bundle $\mathfrak g=VG|_I$ on $I$ is the \emph{(differentiable ) algebroid}
associated to the groupoid $G$. 
\end{df}

Given sections $\xi$, $\eta$ of $s:\mathfrak g\rightarrow I$ defined on an open set $U$ of  $I$, it is well defined the Lie bracket $\zerobra{\,\,}{\,}$ on local sections of $ \mathfrak g$, given by
\begin{equation}\label{Liebra}
\zerobra{\xi}{\eta}=\left[\overline\xi,\overline\eta\right]|_I.
\end{equation}
\begin{po}
If $f$ is a real function on $U$, $\xi$, $\eta$ sections of $\mathfrak g$ on $U$, then
$$
\zerobra{f\xi}{\eta}=f\zerobra{\xi}{\eta}-(t_*\eta)(f)\xi.
$$
\end{po}
\Pf As $\overline {f\xi}=(f\circ t)\overline\xi$, it follows
$$
\zerobra{f\xi}{\eta}=\left[(f\circ t)\overline\xi,\overline\eta\right]|_I=\left((f\circ t)\left[\overline\xi,\overline\eta\right]-\overline\eta(f \circ t)\overline\xi\right)|_I=f\zerobra{\xi}{\eta}-(t_*\eta)(f)\xi.
$$\EPf

If $\Gamma(\mathfrak{g})$ denotes the sheaf of germs  of local sections of $\mathfrak{g}$, then $\Gamma(\mathfrak{g})$ is a Lie algebra sheaf, with the Lie bracket $\zerobra{\,\,}{\,\,}$.

\begin{po}\label{difinverse}
If $\iota$  denotes the inverse in $G$  and if $\xi\in \mathfrak g_x$, then 
$$
\iota_*(\xi)=t_*(\xi)-\xi.
$$
\end{po}
\Pf Consider $X_\epsilon$ a curve in $G_x$ such that $\frac{d}{d\epsilon}X_\epsilon |_{\epsilon=0}=\xi$. Denote by $f(\epsilon)=t(X_\epsilon)$. The composition $X_\epsilon.(X_\epsilon)^{-1}=I_{f(\epsilon)}$ so
$$
\frac{d}{d\epsilon}(X_\epsilon.X_\epsilon^{-1} )|_{\epsilon=0}=\frac{d}{d\epsilon}I_{f(\epsilon)} |_{\epsilon=0}.
$$
Observe that  $\xi$ is $s$-vertical and $\iota_*(\xi)$ is $t$-vertical.
Therefore 
$$
\frac{d}{d\epsilon}X_\epsilon|_{\epsilon=0}.I_x+I_x.\frac{d}{d\epsilon}X_\epsilon^{-1}|_{\epsilon=0}=\frac{d}{d\epsilon}I_{f(\epsilon)} |_{\epsilon=0}
$$
and we get
$$
\xi+\iota_*(\xi)=t_*(\xi).
$$
\EPf

%XXXXXXXXXXXXXXXXXXX
%XXXXXXXXXXXXXXXXX
%XXXXXXXXXXXXXXXX

\subsection{Right invariant diffeomorphisms of $G$}

\begin{df}
A diffeomorphism $F:G\rightarrow G$ is right invariant if $F(XY)=F(X)Y$ for every $(X,Y)\in(s\times t)^{-1}\Delta$.
\end{df}
Every right invariant diffeomorphism is defined by a bisection $\sigma:I\rightarrow G$ defined as $\sigma(x)=F(I_x)$. In fact
$$
F(X)=F(I_{t(X)}X)=F(I_{t(X)})X=\sigma(t(X))X.
$$
From this formula it follows that $F$ takes $G(.,t(X))$ on $G(.,f(t(X)))$ where $f=tF$. As $F$ is a diffeomorphism, the function $f$ must be a diffeomorphism and  $\sigma$ must be a bisection.

The composition of right invariant diffeomorphisms is a right invariant diffeomorphism. In fact, if $\tilde F$ is a right invariant diffeomorphism such that $\tilde F(Y)=\tilde\sigma(t(Y))Y$, then
$$\tilde F\circ F(X)=\tilde F(F(X))=\tilde\sigma(t(F(X)))F(X)=\tilde\sigma(t(\sigma(t(X))))\sigma(t(X))X=\tilde\sigma\circ\sigma(t(X))X.
$$
\begin{po}\label{f41}
If $v$ is a vector field tangent to $I$, and $\sigma\in\mathcal G$ is such that $t\sigma=f$, then 
$$
j^1\sigma.v.j^1\sigma^{-1}=f_*(v).
$$
\end{po}
\Pf 
Posing $v=\frac{d}{du}x_u|_{u=0}$, we obtain
$$
j^1\sigma.v.j^1\sigma^{-1}=\frac{d}{du}(\sigma(x_u).\sigma^{-1}(f(x_u)))|_{u=0}=\frac{d}{du}f(x_u)|_{u=0}=f_*v,
$$\EPf

%%%XXXXXXXXXXXXXXXXXXXXX

%XXXXXXXXXXXXXXXXXXXXX

%XXXXXXXXXXXXXXXXXXXXX

\subsection{Groupoids and algebroids of jets}
Let be $G$ a differentiable groupoid and $G_k$ the manifold of $k$-jets of local bisections  of $G$. This manifold has a natural structure of Lie groupoid given by composition of jets
$$j^k_{f(x)}H.j^k_xF=j^k_x(HF),$$
and inversion
$$(j^k_xF)^{-1}=j^k_{f(x)}F^{-1},$$
where  $F:U\rightarrow G$,  $H:V\rightarrow G$ are local bisections of $G$, $f=tF$, $V=f(U)$ and $x\in U$. The groupoid $G_k$ has a natural submanifold of identities $I_k=j^k\mbox{I}$, where $\mbox{I}$ is the identity section of $G$. We have a natural identification of $I$ with $I_k$, given by $I(x)\mapsto I_k(x)$. Therefore we can think of $I$ as a submanifold of $G_k$. There are two submersions  $s,t:G_k\rightarrow I$, the canonical projections \emph{source}, $s(j^k_xF)=x$, and \emph{target}, $t(j^k_xF)=f(x)$.  We also consider $s$ and $t$ with values in $I$, by the above identification of $I_k$ with $I$.

There are natural projections $\pi^k_l=\pi_l:G_k\rightarrow G_l$, for $l\geq 0$, defined by $\pi_l(j^k_xF)=j^l_xF$. Observe that $G_0=G$. The projections $\pi^k_l$ commute with the operations of composition and inversion in $G_k$. A bisection $F_k$ of $G_k$ is \emph{holonomic} if there exist a bisection $F$ of $G$ such that
$F_k=j^kF$. Therefore, if $F_k$ is holonomic, we have $F_k=j^k(\pi_0F_k)$.

\begin{df}
The vector bundle $\mathfrak g_k=VG_k|_I$ on $I$ is the \emph{(differentiable) algebroid}
associated to the groupoid $G_k$. 
\end{df}

 The vector bundle $J^k\mathfrak g$ is a vector bundle on $I$, and we also denote by $\pi:J^k\mathfrak g\rightarrow  I$ the map $\pi(j^k_x\theta)=x$, where  $\theta:U\subset I\rightarrow \mathfrak g$ is a local section. If $f_t$ is the 1-parameter group of local diffeomorphisms of $G$ such  that $\frac{d}{dt}f_t|_{t=0}=\theta$, then we get, for $x\in U$, $$\frac{d}{dt}j^k_xf_t|_{t=0}=j^k_x\theta.$$
  This means we have a natural identification 
  $$\mathfrak g_k=J^k\mathfrak g.$$
  
  \begin{po}\label{bracketk}
The bracket $\zerobra{\,\,}{\,\,}_k$ on $\Gamma({\mathfrak g}_k)$ is determined by: 
\begin{enumerate}
	\item[(i)] $\zerobra{j^k\xi}{j^k\eta}_k=j^k[\xi,\eta],\ \  \xi,\eta\in \Gamma(\mathfrak g)$
	\item[(ii)] $\zerobra{\xi_k}{f\eta_k}_k=f\zerobra{\xi_k}{\eta_k}_k+(t_*\xi_k)(f)\eta_k$, 
\end{enumerate}
where $\xi_k,\eta_k\in \Gamma(\mathfrak g_k)$, $f$ a real function on $I$.
\end{po}
 
 Let be $T=TI$ the tangent bundle of $I$, and  $\mathcal T$ the sheaf of germs of local sections of $T$. 
   Therefore, as $TG_k|_I=TI\oplus VG_k|_I$,  
   \begin{equation}\label{tqki}
   TG_k|_I\cong T\oplus \mathfrak g_k,
   \end{equation}
 and if we denote by
 $$\check{J}^k\mathfrak g=T\oplus J^k\mathfrak g,$$
then $TQ^k|_I\cong \check{J}^k\mathfrak g$. Observe that $\check{J}^k\mathfrak g $ is a vector bundle on $I$. The restriction of $t_*:TG_k\rightarrow T$ to $TG_k|_I$, and the isomorphism $TG_k|_I\cong \check{J}^k\mathfrak g$ defines the map, that we denote again by $t_*$, 
 \begin{equation}\label{beta}
 \begin{array}{rcl}t_*:\check{J}^kT&\rightarrow &T\\
 v+j_x^k\theta\in (\check{J}^kT)_x&\mapsto &v+t_*\theta(x)
 \end{array}.
 \end{equation}

We denote by the same symbols as above the projections $\pi^k_l=\pi_l:J^k\mathfrak g\rightarrow J^l\mathfrak g$, $l\geq 0$, defined by $\pi_l(j^k_x\theta)=j^l_x\theta$. If $\xi_k$ is a point or a section of $J^k\mathfrak g$, we denote by $\xi_l=\pi^k_l(\xi_k)$. The vector bundle $J^0\mathfrak g$ is isomorphic to $\mathfrak g$.
 
We have the canonical inclusions
\begin{equation*}
\begin{array}{rcl}\label{lkl}
\lambda^l:J^{k+l}\mathfrak g&\rightarrow &J^lJ^k\mathfrak g\\
j^{k+l}_x\theta&\mapsto&j^l_x(j^k\theta)
\end{array}.
\end{equation*}
for $\theta\in\Gamma(\mathfrak g)$. 

Analogously to definition of bisections of $G_k$, a section $\xi_k$ of $J^k\mathfrak g$ is \emph{holonomic} if there exist $\xi\in\Gamma(\mathfrak g)$ such that
$\xi_k=j^k\xi$. Therefore, if $\xi_k$ is holonomic, we have $\xi_k=j^k(\pi_0\xi_k)$.

We denote by $G_{k,l}$ the groupoid
$$
G_{k,l}=\{j^l_xF:x\in M, F \mbox{ a bisection of } G_k \}
$$
and by $\pi^{k,l}_{k',l'}=\pi_{k',l'}:G_{k,l}\rightarrow G_{k',l'}$ the natural projections $\pi^{k,l}_{k',l'}(j^l_xF)=j^{l'}_x(\pi_{k'}F)$.

The Lie algebroid associated to $G_{k,l}$ is $J^l(J^k\mathfrak g)$ which we denote also by $J^{k,l}\mathfrak g$.
%XXXXXXXXXXXXXXXXXXXXXXXXXXXXXXXXXXXXXXXXXX

\subsection{The affine structures on $J^1G$ and $J^1\mathfrak g$}
Given two sections $\sigma,\eta$ of $\mathfrak g$ such that $\sigma(x)=\eta(x)$ and $v\in T_x$ we have $\sigma_*v-\eta_*v\in V_{\eta(x)}\mathfrak g.$ As $\mathfrak g$ is a vector bundle, $V\mathfrak g\cong \mathfrak g$, therefore $\sigma_*v-\eta_*v\in\mathfrak g_x.$ This means $j^1_x\sigma-j^1_x\eta\in T^*\otimes \mathfrak g$. Inversely, given $j^1_x\sigma$ and $u\in T^*\otimes \mathfrak g$ there exists $\eta$ section of $\mathfrak g$ such that $j^1_x\sigma-j^1_x\eta=u$. As $J^1\mathfrak g$ has the canonical $0$-section we get the exact sequence
\begin{equation}\label{exact1}
0\rightarrow T^*\otimes \mathfrak g\rightarrow J^1\mathfrak g\rightarrow \mathfrak g\rightarrow 0.
\end{equation}

The same argument  applies to $J^1G$. In this case, if $\sigma,\eta$ are sections of $G$ such that $\sigma(x)=\eta(x)$ and $v\in T_x$ we have $\sigma_*v-\eta_*v\in V_{\eta(x)}G$, where 
$$
VG=\{v\in TG:t_*v=0\}.
$$ 
This means $j^1_x\sigma-j^1_x\eta\in T^*_x\otimes V_{\eta(x)}G$. Inversely, given $j^1_x\sigma$ and $u\in T^*_x\otimes V_xG$ there exists $\eta$ section of $G$ such that $j^1_x\sigma-j^1_x\eta=u$. A special case occurs when $\eta$ is the identity section $I$. As $VG|_I=\mathfrak g$, in this case we get $j^1_x\sigma-j^1_xI\in T^*_x\otimes \mathfrak g_x$. Inversely, given  $u\in T^*_x\otimes \mathfrak g_x$ there exists $\sigma$ section of $G$ such that $j^1_x\sigma-j^1_xI=u$. 

It is important to characterize the set of $u\in T^*_x\otimes \mathfrak g_x$ such that $j^1_x\sigma-j^1_xI=u\in T_x^*\otimes \mathfrak g_x$ is such that $\sigma$ is a bisection of $G$:
\begin{po}\label{image}
The set of $u\in T^*_x\otimes \mathfrak g_x$ such that $j^1_x\sigma =j^1_xI+u$ is the jet of a bisection $\sigma$ of $G$ is characterized by the application
$$
v\in T_x\mapsto v+t_*(i(v)u)\in T_x
$$
to be invertible.
\end{po}
\Pf The proof is an easy consequence of  $\sigma$ being a bisection if and only if $t\sigma$ is a diffeomorphism. \EPf

As $J^k\mathfrak g\subset J^1J^{k-1}\mathfrak g$ we obtain as a particular case of (\ref{exact1})  the exact sequence
\begin{equation*}\label{exact2}
0\rightarrow S^kT^*\otimes \mathfrak g\rightarrow J^k\mathfrak g\rightarrow J^{k-1}\mathfrak g\rightarrow 0.
\end{equation*}
For a proof see \cite{N}.

%XXXXXXXXXXXXXX
\subsection{The linear Spencer operator}
If $\theta:U\subset I\rightarrow J^k\mathfrak g$ is a section, let be $\xi=j^1_x\theta\in J^1_xJ^k\mathfrak g$, $x\in U$. Then $\xi$ can be identified to a linear application
\begin{equation*}
\begin{array}{rcl}
	\xi_*:T_xI&\rightarrow &T_{\theta(x)}J^k\mathfrak g\\
	v&\mapsto &\theta_*(v)
	\end{array}.
\end{equation*}
If $\eta\in J^1_xJ^k\mathfrak g$ is given by $\eta=j^1_x\mu$, with $\mu(x)=\theta(x)$,   then $(\pi)_*(\eta_*-\xi_*)v=0$, where we remember that $\pi:J^k\mathfrak g\rightarrow I$ is defined by $\pi(j^k_x\theta )=x$.
So $\eta_*-\xi_*\in T_x^*I\otimes V_{\pi^1_0(\xi)}J^k\mathfrak g$, where $VJ^k\mathfrak g=\ker\pi_{*}$.   But $J^k\mathfrak g$ is a vector bundle, then $V_{\pi^1_0\xi}J^k\mathfrak g\cong J^k_x\mathfrak g$, so $\eta_*-\xi_*\in T_x^*I\otimes J^k_x\mathfrak g$. The sequence
\begin{equation}\label{kernelpi}
0\rightarrow T^*I\otimes J^k\mathfrak g\rightarrow J^1J^k\mathfrak g\stackrel{\pi^1_0}{\rightarrow}J^k\mathfrak g\rightarrow 0
\end{equation}
obtained in this way is exact.

The linear operator  $D$  defined by
\begin{equation}\label{D}
\begin{array}{rcl}
	D:\Gamma(J^k\mathfrak g)&\rightarrow &\mathcal{T}^*\otimes \Gamma(J^{k-1}\mathfrak g)\\
	\xi_k&\mapsto &D\xi_k=j^1\xi_{k-1}-\lambda^1(\xi_k),
	\end{array}
\end{equation}
is the \emph{linear Spencer operator}. We remember that $\xi_{k-1}=\pi^k_{k-1}\xi_k$ and 
\begin{equation*}
\begin{array}{rcl}
	\lambda^1:J^k\mathfrak g&\rightarrow &J^1J^{k-1}\mathfrak g\\
	j^k_x\xi&\mapsto &j^1_x(j^{k-1}\xi)
	\end{array}.
\end{equation*}
The difference in (\ref{D}) is  in $J^1J^{k-1}\mathfrak g$ and is in $T^*\otimes J^{k-1}\mathfrak g$, by (\ref{kernelpi}). 

The operator $D$ is null on a section $\xi_k$ if and only if it is holonomic, i.e., $D\xi_k=0$ if and only if there exists $\theta\in\Gamma(\mathfrak g)$ such that $\xi_k=j^k\theta$. 
\begin{po}\label{propriedadesD}The operator $D$ is characterized by
\begin{enumerate}
	\item[(i)] $D\circ j^k=0$\
	\item[(ii)] $D(f\xi_k)=df\otimes \xi_{k-1}+fD\xi_k$, 
\end{enumerate}
with $\xi_k\in \Gamma(J^k\mathfrak g)$, $\xi_{k-1}=\pi_{k-1}\xi_k$ and $f$ a real function on $I$. 
\end{po}
For a proof, see \cite{KS}.

The operator $D$ extends to
\begin{equation*}\label{extD}
\begin{array}{rcl}
	D:\wedge^l\mathcal{T}^*\otimes \Gamma(J^k\mathfrak g)&\rightarrow &\wedge^{l+1}\mathcal{T}^*\otimes \Gamma(J^{k-1}\mathfrak g)\\
	\omega\otimes\xi_k&\mapsto &D(\omega\otimes\xi_k)=d\omega\otimes \xi_{k-1}+(-1)^l\omega\wedge D\xi_k.
	\end{array}
\end{equation*}
%XXXXXXXXXXXXXX

\section{The calculus on the diagonal}\label{calcdiag}

Following \cite{Ma1}, \cite{Ma2}, \cite{KS} and \cite{V}, we will relate $\check{J}^k\mathfrak g$ with vector fields along the diagonal of $I\times G$ and actions of bisections in $G_k$ with diffeomorphisms of $I\times G$ which are right invariants,  $\rho_1$ projectables and preserve $B(G)=\{(tX,X):X\in G\}$. 

We denote the \emph{diagonal} of $I\times G$ by $\Delta=\{(x,x)\in I\times G|x\in I\}$, and by $\rho_1:I\times G\rightarrow I$ and $\rho_2:I\times G\rightarrow G$ the first and second projections, respectively. The restrictions $\rho_1|_\Delta$, $\rho_2|_\Delta$ and $t\circ\rho_2|_\Delta$ are diffeomorphisms of $\Delta$ on $I$. A sheaf on $I$ will be identified to its inverse image by $\rho_1|_\Delta$. For example, if $\mathcal{O}_I$ denotes the sheaf of germs  of real functions on $I$, then we will write $\mathcal{O}_I$ on $\Delta$ instead $(\rho_1|_{\Delta})^{-1}\mathcal{O}_I$. Therefore, a $f\in \mathcal{O}_I$ will be considered in $\mathcal{O}_{\Delta}$ or in $\mathcal{O}_{I\times G}$ through the map $f\mapsto f\circ \rho_1$. 

The transposition in $I\times I$ is denoted by
\begin{equation}
\begin{array}{rcl}\label{trans}
\epsilon:I\times I&\rightarrow & I\times I\\
(x,y)&\rightarrow&(y,x).
\end{array}
\end{equation}

The right action of $G$ on $G$ extends  to $I\times G$ by 
$$
\begin{array}{rcl}
(I\times G)\times G&\rightarrow & I\times G\\
((a,Y),X)&\rightarrow&(a,YX)
\end{array}
$$
where $(X,Y)\in (s\times t) ^{-1}(\Delta)$.
A vector field $\xi$ on $I\times G$ is \emph{right invariant} if is tangent to the submanifolds $I\times s^{-1}(x)$ for every $x\in I$ and $\xi(a,XY)=\xi(a,X)Y$. A right invariant vector field on $I\times G$ is defined by its restriction to $I\times I$ since that $\xi(a,X)=\xi(a,t(X))X$.

%XXXXXXXXX

%XXXXXXXX

\subsection{Brackets in $\check{J}^k\mathfrak g$ }
We denote by $\mathcal{T}(I\times G)$ the sheaf of germs of local sections of $T(I\times G)\rightarrow I\times G$; by $\mathcal{R}$ the sub sheaf in Lie algebras of $\mathcal{T}({I}\times {G})$ of right invariant vector fields whose elements are $\rho_1$-projectables; by $\mathcal{H_R}$ the sub sheaf in Lie algebras of $\mathcal{R}$ that projects on $0$ by $\rho_2$, i.e. $\mathcal{H_R}=(\rho_2)_*^{-1}(0)\cap\mathcal{R}$; and by $\mathcal{V_R}$ the sub sheaf in Lie algebras defined by $\mathcal{V_R}=(\rho_1)_*^{-1}(0)\cap\mathcal{R}$. Clearly, $$\mathcal{R}=\mathcal{H_R}\oplus \mathcal{V_R},$$ and $$[\mathcal{H_R},\mathcal{V_R}]\subset\mathcal{V_R}.$$
 Then 
$$(\rho_1)_*:\mathcal{H_R}\widetilde{\longrightarrow} 
\mathcal {T}$$ 
is an isomorphism, so we identify $\mathcal{H_R}$ naturally with $\mathcal{T}$ by this isomorphism, and utilize both notations indistinctly.
\begin{po}\label{bili}
The Lie bracket in $\mathcal{R}$ satisfies: 
\begin{equation*}
[v+\xi,f(w+\eta)]=v(f)(w+\eta)+f[v+\xi,w+\eta],
\end{equation*}
with $f\in\mathcal{O}_I$, $v,w\in \mathcal{H_R}$, $\xi,\eta\in\mathcal{V_R}$.
In particular, the Lie bracket in $\mathcal{V_R}$ is $\mathcal{O}_I$-bilinear.\end{po} 
\Pf Let be $f\in \mathcal{O}_I$, $\xi,\eta\in\mathcal{V_R}$. Then 
$$[v+\xi,(f\circ\rho_1)(w+\eta)]=(v+\xi)(f\circ\rho_1)(w+\eta)+(f\circ\rho_1)[v+\xi,w+\eta].$$
As $f\circ\rho_1$ is constant on the sub manifolds $\{x\}\times G$ and
$\xi$ is tangent to them, we obtain $\xi(f\circ\rho_1)=0$, and the proposition follows.
\EPf

We know that a right invariant vector field on $G$ is defined by its restriction to $I$, therefore identifies to an element of $\Gamma(\mathfrak g)$.
A vector field in $\mathcal{V_R}$ is given by a family of sections of $\Gamma(\mathfrak g)$ parameterized by an open set of $I$. 
%Let's denote by 	$\mathcal{R}|_{\Delta}$ the subset of $	\mathcal{R}$ whose elements are sections defined on a neighborhood of $\Delta$. 
Therefore there exists a surjective morphism
\begin{equation*}\label{R}
\begin{array}{rcl}
	\Upsilon_k:\mathcal{R}&\rightarrow &\mathcal{T}\oplus \Gamma(J^k\mathfrak g)\\
	v+\xi&\mapsto &v+\xi_k,
	\end{array}
\end{equation*}
where $v\in\mathcal{H_R}$, $\xi\in\mathcal{V_R}$, and $$\xi_k(x)=j^k_{(x,x)}(\xi|_{\{x\}\times I}).$$ 
The kernel of this morphism is the sub sheaf $\mathcal{V_R}^{k+1}$ of $\mathcal{V_R}$ constituted by vector fields that are null on $\Delta$ up to order $k$. Therefore $\mathcal{R}/\mathcal{V_R}^{k+1}$ is null outside of $\Delta$, and can be considered as the sections of a vector bundle on $\Delta$, and this vector bundle is isomorphic to the vector bundle on $I$, ${T}\oplus J^k\mathfrak g$. Observe that the sections considered in the quotient are sections on open sets of $I$.
As $$[\mathcal{R},\mathcal{V_R}^{k+1}]\subset \mathcal{V_R}^k,$$
the bracket 
on $\mathcal{R}$ induces a bilinear antisymmetric map which we call the  \emph{first bracket} of order $k$,
\begin{equation}\label{jka}
\firstbra{\,\,}{\,\,}_k=(\mathcal{T}\oplus \Gamma(J^k\mathfrak g))\times (\mathcal{T}\oplus \Gamma(J^k\mathfrak g))\rightarrow \mathcal{T}\oplus \Gamma(J^{k-1}\mathfrak g)
\end{equation}
defined by 
$$
\firstbra{v+\xi_k}{w+\eta_k}_k=\Upsilon_{k-1}([v+\xi,w+\eta]),
$$
where $\Upsilon_k(\xi)=\xi_k$ and $\Upsilon_k(\eta)=\eta_k$. 

It follows from proposition \ref{bili} that $\firstbra{\,\,}{\,\,}_k$ satisfies:
\begin{equation}\label{noombi}
\firstbra{v+\xi_k}{f(w+\eta_k)}_k=v(f)(w+\eta_{k-1})+f\firstbra{v+\xi_k}{w+\eta_k}_k,
\end{equation}

\begin{eqnarray*}\label{jacobipointbracket}
\firstbra{\firstbra{v+\xi_k}{w+\eta_k}_k}{z+\theta_{k-1}}_{k-1}&+&
\firstbra{\firstbra{w+\eta_k}{z+\theta_k}_k}{v+\xi_{k-1}}_{k-1}\nonumber\\
&+&\firstbra{\firstbra{z+\theta_k}{v+\xi_k}_k}{w+\eta_{k-1}}_{k-1}=0,
\end{eqnarray*}

for $v,w,z\in\mathcal{T}$, $\xi_k,\eta_k,\theta_k\in \Gamma(J^k\mathfrak g)$, $f\in\mathcal{O}_I$. In particular, the first bracket
is $\mathcal{O}_I$-bilinear on $\Gamma(J^k\mathfrak g)$. Also
$$\firstbra{\Gamma(J^0\mathfrak g)}{\Gamma(J^0\mathfrak g)}_0=0.$$

The following proposition relates $\firstbra{\,\,}{\,\,}_k$ with the bracket in $\mathcal{T}$ and the linear Spencer operator $D$ in $\Gamma(J^k\mathfrak g)$.
%There exists also a bilinear antisymmetric map in $J^kT$ with values in $J^{k-1}T$ defined by
%\begin{equation}\label{pointbracket}
%\begin{array}{rcl}
%	\firstbra{\,\,}{\,}:J^kT\times J^kT&\rightarrow & J^{k-1}T\\
%	(j^k_x\xi,j^k_x\eta)&\mapsto &\firstbra{j^k_x\xi}{j^k_x\eta}_k=j^{k-1}_x[\xi,\eta].
%	\end{array}
%\end{equation}
%Clearly, if $\xi_k,\eta_k,\theta_k\in J^k\T$, $f,g$ real functions on $M$, then:
%\begin{equation}\label{jacobipointbracket}
%\firstbra{\firstbra{\xi_k}{\eta_k}_k}{\theta_{k-1}}_{k-1}+
%\firstbra{\firstbra{\eta_k}{\theta_k}_k}{\xi_{k-1}}_{k-1}+
%\firstbra{\firstbra{\theta_k}{\xi_k}_k}{\eta_{k-1}}_{k-1}=0,
%\end{equation}
%\begin{equation}\label{flinearpointbracket}
%\firstbra{f\xi_k}{g\eta_k}_k =fg\firstbra{\xi_k}{\eta_k}_k.
%\end{equation}
\begin{po}\label{colchetek}
Let be $v,w\in\mathcal{T}$, $\theta,\mu\in\Gamma(\mathfrak g)$, $\xi_k,\eta_k\in \Gamma(J^k\mathfrak g)$ and $f\in\mathcal{O}_I$. Then:
\begin{enumerate}
\item [(i)] $\firstbra{v}{w}_k=[v,w]$, where the  bracket at right is the bracket in $\mathcal{T}$;
\item[(ii)] $\firstbra{v}{\xi_k}_k=i(v)D\xi_k$;
\item[(iii)] $\firstbra{j^k\theta}{j^k\mu}_k=j^{k-1}[\theta,\mu]$, where the bracket at right is the bracket in $\Gamma(\mathfrak g)$.
\end{enumerate}
\end{po}
\Pf (i) This follows from the identification of $\mathcal{T}$ with $\mathcal{H_R}$.

(ii) First of all, if $\theta\in\Gamma(\mathfrak g)$, let be $\Theta\in\mathcal{V_R}$ defined by $\Theta(x,Y)=\theta(t(Y))Y$. Then $\Upsilon_k(\Theta)=j^k\theta$. If $v\in\mathcal{H_R}$, then $v$ and $\Theta$ are both $\rho_1$ and $\rho_2$ projectables,  $(\rho_1)_*(\Theta)=0$ and $(\rho_2)_*(v)=0$, so we get $[v,\Theta]=0$. Consequently 
\begin{equation}\label{bvt}
\firstbra{v}{j^k\theta}_k=\Upsilon_{k-1}([v,\Theta])=0.
\end{equation}
Also by (\ref{noombi}) we have
\begin{equation}\label{vfx}
\firstbra{v}{f\xi_k}_k=v(f)\xi_{k-1}+f\firstbra{v}{\xi_k}_k.
\end{equation}
As (\ref{bvt}) and (\ref{vfx}) determine $D$ (cf. proposition \ref{propriedadesD}), we get (ii).

(iii) Given $\theta,\mu\in\Gamma(\mathfrak g)$, we define $\Theta,H\in\mathcal{V_R}$ as in (ii), $\Theta(x,Y)=\theta(t(Y))Y$ and $H(x,Y)=\mu(t(Y))Y$. Therefore
$$
\firstbra{j^k\theta}{j^k\mu}_k=\firstbra{\Upsilon_k\Theta}{\Upsilon_kH}_k=\Upsilon_{k-1}([\Theta,H])=j^{k-1}[\theta,\mu].
$$
%Choosing a set of $\theta$'s such that the $\Theta_k$'s are a basis of $J^k\mathcal{A}$ as an $\mathcal{O}_I$-module, and considering the $\mathcal{O}_I$-bilinearity of $[\,\, ,\,\,]_k$
%and $[\,\, ,\,\,]_k$, (iii) follows.
\EPf

Let be $\tilde{\mathcal{V}}_\mathcal{R}$ the sub sheaf in Lie algebras of $\mathcal{R}$ such that $\tilde\xi\in\tilde{\mathcal{V}}_\mathcal{R}$ if and only if $\tilde\xi$ is tangent to the submanifold $B(G)$ of $I\times G$ image of $G$ by the injective function 
$$\begin{array}{rcl}
%t\times \mbox{id}
B:G&\rightarrow &I\times G\\
X&\mapsto &(t(X),X),
\end{array}
$$
i.e, $B(G)=\{(t(X),X)\in I\times G|X\in G\}$.
If $\tilde\xi=\xi_H+\xi\in\tilde{\mathcal{V}}_\mathcal{R}$, $\xi_H\in\mathcal{H}_\mathcal{R}$, $\xi\in\mathcal{V}_\mathcal{R}$, then  
$$\xi_H(t(X),X)+\xi(t(X),X)=\tilde\xi(t(X),X).$$
If $X(u)$, $u\in (-a,a)$, is a curve such that $\frac{d}{du}X(u)|_{u=0}=\xi$, then 
$$\tilde\xi= \frac{d}{du}(t(X(u)),X(u))|_{u=0}=\epsilon_*(t_*\xi)+\xi.$$
Therefore $\xi_H=\epsilon_*(t_*\xi)$. Remember that $\epsilon$ is the transposition (\ref{trans}).

%As $$t_*(\xi(t(X),X))=t_*(\xi(t(X),t(X)).X)=t_*(\xi(t(X),t(X))),$$
%it follows that for $x\in I$,
 %\begin{equation}\label{xiH}
% (\rho_1)_*(\xi_H(x,x))=(t\circ \rho_2)_*(\xi(x,x)).
 %\end{equation}
 %where $\rho_1,\rho_2:M\times M\rightarrow \Delta$. 
 Consequently, if $\xi_k=\Upsilon_k(\xi)$, then $\xi_H=\epsilon_*t_*(\xi_k)$, where we remember that  $t_*:J^k\mathfrak g\rightarrow T$ is defined in (\ref{beta}). 
 %We can, too, write $\xi_H=(\rho_2)_*(\xi_k)$, since that $(\rho_2)_*\mathcal{V_R}^{k+1}|_{\Delta}=0$. 
 From now on, $\xi_H$ denotes the horizontal component of $\tilde\xi\in\tilde{\mathcal{V}}_\mathcal{R}$, so $\tilde\xi=\xi_H+\xi$, with $\xi_H\in\mathcal{H}$ and $\xi\in\mathcal{V}$.
 
 We denote by $\Gamma(\tilde J^k\mathfrak g)$ the sub sheaf of $\mathcal{T}\oplus \Gamma(J^k\mathfrak g)$, whose elements are 
 $$\tilde\xi_k=\xi_H+\xi_k,$$
  where $\xi_H=\epsilon_*t_*(\xi_k)$. 
  %or $\xi_H=(\rho_2)_*(\xi_k)$
 If $\eta$ is a vector field in $\mathcal{V_R}^{k+1}$ and if $(x,x)\in\Delta$, then $\eta(x)=0$ up to order $k$. It follows from (\ref{actionvert}) that if $t(X)=x$ then $\eta_X=\eta_x.X$ is null of order $k$ at $X$. Therefore $\eta$ is null of order $k$ on $B(G)$ and 
 $\mathcal{V_R}^{k+1}\subset\tilde{\mathcal{V}}_\mathcal{R}$. Furthermore $\Gamma(\tilde J^k\mathfrak g)$ identifies with $\tilde{\mathcal{V}}_\mathcal{R}/\mathcal{V_R}^{k+1}$ since that 
 \begin{equation}\label{VtV}
 [\tilde{\mathcal{V}}_\mathcal{R},\mathcal{V_R}^{k+1}]\subset \mathcal{V_R}^{k+1},
 \end{equation}
  because the vector fields in $\tilde{\mathcal{V}}_\mathcal{R}$ are tangents to $B(G)$. It follows that the bracket in $\tilde{\mathcal{V}}_\mathcal{R}$ defines a bilinear antisymmetric map,  called the \emph{second bracket}, by
 \begin{equation*}\label{bta} 
 \begin{array}{rcl}
	\secondbra{\,\,}{\,\,}_k:\Gamma(\tilde J^k\mathfrak g)\times \Gamma(\tilde J^k\mathfrak g)&\rightarrow &\Gamma(\tilde J^k\mathfrak g) \\
	(\xi_H+\xi_k,\eta_H+\eta_k)&\mapsto &\Upsilon_k([\tilde\xi,\tilde\eta]),
	\end{array}
 \end{equation*}
 where $\Upsilon_k(\tilde\xi)=\xi_H+\xi_k$ and $\Upsilon_k(\tilde\eta)=\eta_H+\eta_k$.
Unlike the first bracket (\ref{jka}), we do not loose one order doing the bracket in $\tilde J^k\mathfrak g$. 
%To be more precise, if $\xi_{k+1}$ and $\eta_{k+1}$ are relevements of $\xi_k$ and $\eta_k$, then
The second bracket $\secondbra{\,\,}{\,\,}_k$ is a Lie bracket on $\tilde J^k\mathfrak g$. The  proposition \ref{tilde} below relates it with the  bracket $\zerobra{\,\,}{\,\,}_k$, defined in (\ref{Liebra}).

The projection 
$$\begin{array}{rcl}
\nu:\mathcal{H_R}\oplus\mathcal{V_R}&\rightarrow &\mathcal{V_R}\\
v+\xi&\mapsto &\xi
\end{array}
$$
quotients to
\begin{equation*}%\label{nu}
\begin{array}{rcl}	\nu_k:\mathcal{T}\oplus \Gamma( J^k\mathfrak g)&\rightarrow &\Gamma( J^k\mathfrak g) \\
v+{\xi}_k	&\mapsto  &\xi_k.
	\end{array}
\end{equation*}
and the restriction $\nu_k:\tilde J^k\mathfrak g\rightarrow J^k\mathfrak g$ is an isomorphism of vector bundles.
\begin{po}\label{tilde}
If $\tilde\xi_k,\tilde\eta_k\in \Gamma(\tilde J^k\mathfrak g)$ then 
$$\zerobra{\xi_k}{\eta_k}_k=\nu_k(\secondbra{\tilde\xi_k}{\tilde\eta_k}_k),$$
where $\xi_k=\nu_k(\tilde\xi_k)$, $\eta_k=\nu_k(\tilde\eta_k)$.
\end{po}
\Pf We will verify properties (i) and (ii) of Proposition \ref{bracketk}. If $\theta,\mu\in\mathfrak g $, let be $\Theta,H\in\mathcal{V_R}$ as in the proof of proposition \ref{colchetek}. Then:
\begin{enumerate}
\item[(i)] 
$\begin{array}{rcl}
\nu_k(\secondbra{\epsilon_*t_*\theta+j^k\theta}{\epsilon_*t_*\mu+j^k\mu}_k)&=&\nu_k(\Upsilon_k([\epsilon_*t_*\theta+\Theta,\epsilon_*t_*\mu+H]))\vspace{.2cm}\\
&=&\Upsilon_k(\nu(\epsilon_*t_*[\theta,\mu]+[\Theta,H]))=j^k([\theta,\mu])=\zerobra{j^k\theta}{j^k\mu}_k.
\end{array}$

\item[(ii)] 
$\begin{array}{rcl}	\nu_k(\secondbra{\tilde\xi_k}{f\tilde\eta_k}_k)&=&\nu_k(f\secondbra{\tilde\xi_k}{\tilde\eta_k}_k+\xi_H(f)\tilde\eta_k)\vspace{.2cm} \\
	&= &f\nu_k(\secondbra{\tilde\xi_k}{\tilde\eta_k}_k)+(\epsilon_*t_*\xi_k)(f)\eta_k.
	%&=&[[\xi_k,f\eta_k]]_k.
	\end{array}
$
\end{enumerate}
\EPf
\begin{co}
If $\tilde\xi_k,\tilde\eta_k\in \Gamma(\tilde J^k\mathfrak g)$, then 
$$\nu_k\secondbra{\tilde\xi_k}{\tilde\eta_k}_k=\zerobra{\xi_k}{\eta_k}_k=i(\xi_H)D\eta_{k+1}-i(\eta_H)D\xi_{k+1}+\firstbra{\xi_{k+1}}{\eta_{k+1}}_{k+1},$$
where $\xi_H=\epsilon_*t_*\xi_k,\,\eta_H=\epsilon_*t_*\eta_k\in\mathcal{T}$ and $\xi_{k+1},\eta_{k+1}\in J^{k+1}\mathfrak g$ projects on $\xi_k,\,\eta_k$ respectively.
\end{co}
\Pf It follows from Propositions \ref{colchetek} and \ref{tilde}.\EPf

As a consequence of proposition \ref{tilde} we obtain that
\begin{equation*}\label{isonu}
\begin{array}{rcl}	\nu_k:\tilde J^k\mathfrak g&\rightarrow &J^k\mathfrak g \\
\tilde{\xi}_k	&\mapsto  &\xi_k.
	\end{array}
\end{equation*}
is an isomorphism of Lie algebras sheaves, where the bracket in $\tilde J^k\mathfrak g$ is the second bracket $\secondbra{\,\,}{\,\,}_k$ as defined in (\ref{bta}), and the bracket in $J^k{\mathfrak g}$ is the  bracket $\zerobra{\,\,}{\,\,}_k$ as defined in (\ref{Liebra}).

In a similar way, we obtain from (\ref{VtV}) that we can define the \emph{third bracket} as 
 \begin{equation*}\label{tb} 
 \begin{array}{rcl}
	\thirdbra{\,\,}{\,\,}_k:\Gamma(\tilde J^{k+1}\mathfrak g)\times \Gamma(\check{J}^k{\mathfrak g})&\rightarrow &\Gamma(\check{J}^k{\mathfrak g})\vspace{.15cm} \\
	(\xi_H+\xi_{k+1},v+\eta_k)&\mapsto &\Upsilon_k([\tilde\xi,v+\eta]).
	\end{array}
 \end{equation*}
where $\tilde\xi\in\tilde{\mathcal{V}}_{\mathcal R}$, $v+\eta\in\mathcal{R}$. 

\begin{po}\label{propthird}
The third bracket has the following properties:
\begin{enumerate}
\item[(i)] $$\thirdbra{f\tilde\xi_{k+1}}{g\check\eta_k}_k=f\xi_H(g)\check\eta_k-v(f)g\tilde\xi_{k}+fg\thirdbra{\tilde\xi_{k+1}}{\check\eta_k}_k\,;$$
\item[(ii)] $$\thirdbra{\tilde\xi_{k}}{\firstbra{\check\eta_k}{\check\theta_k}_k}_{k-1}=\firstbra{\thirdbra{\tilde\xi_{k+1}}{\check\eta_k}_k}{\check\theta_k}_k+\firstbra{\check\eta_k}{\thirdbra{\tilde\xi_{k+1}}{\check\theta_k}_{k}\,}_k$$
\item[(iii)]$$\thirdbra{\tilde\xi_{k+1}}{\check\eta_k}_k=\firstbra{\tilde\xi_{k+1}}{\check\eta_{k+1}}_{k+1},$$
\end{enumerate}
where $\tilde\xi_{k+1}=\xi_H+\xi_{k+1}\in\Gamma( \tilde{J}^{k+1}\mathfrak{g})$, $\check\theta_k\in \Gamma( \check{J}^{k}\mathfrak g)$, $\check\eta_{k+1} =v+\eta_{k+1}\in\Gamma( \check{J}^{k+1}\mathfrak{g})$, $\tilde\xi_{k}=\pi_k(\tilde\xi_{k+1})$, $\check\eta_{k} =\pi_k(\check\eta_{k+1})$.
\end{po}
\Pf The proof follows the same lines as the proof of proposition \ref{tilde}.
\EPf
%%%%%%
%Observe that 
%\begin{equation}\label{tildeinv}
%[\mathcal{T},\tilde J^{k+1}\mathcal T]_{k+1}\subset \tilde J^{k}\mathcal T.
%\end{equation}
%In fact, if $v\in \mathcal T$ and $\xi=\sum_{l=1}^rf_lj^{k+1}\sigma_l$, where $f_l$ are  functions on $I$ and $\sigma_l$ are bisections of $\mathfrak g$. Then
%$$
%[v,\tilde\xi]_{k+1}=[v,\sum_{l=1}^rf_lt_*\sigma_l]+[v,\sum_{l=1}^rf_lj^{k+1}\sigma_l]=\sum_{l=1}^rv(f_l)(t_*\sigma_l+j^{k}\sigma_l)=\widetilde{i(v)D\xi}.
%$$

\subsection{Action of bisections on $\check{J}^k\mathfrak g$ }

Let's now verify the relationship that exists between the action of  right invariant diffeomorphisms of $I\times G$,   that are $\rho_1$-projectables and preserve $B(G)$, on $\mathcal{R}$ and actions (\ref{actionl}) and (\ref{actionr}) of $G_1G_k$
on $TG_k$. Let be $\sigma$ a (local) right invariant diffeomorphism of $I\times G$ that is $\rho_1$-projectable. Then
$$\sigma(x,Y)=(f(x),\Phi(x,Y))
%=(f(x),\phi(x,\beta(X)).X)
,$$
where $f\in \mbox{Diff }(I)$, $\Phi(x,Y)=\phi_x(t(Y))Y$ and $\phi_x:I\rightarrow G$ is a bisection for all $x\in I$.

%that is, $\sigma$ is defined by $f\in \mbox{Diff }(I)$, and a function
%$$\begin{array}{rcl}	\phi:M&\rightarrow &\mbox{Diff }( M) \\
%x	&\mapsto  &\phi_x.
%	\end{array},
%$$
%The function 

%$\phi_x(y)=\Phi(x,y)$. 
It follows from $\sigma(B(G))\subset B(G)$  that
$$
\sigma(t(X),X)=(f(t(X)),\Phi(t(X),X))=(f(t(X)),\Phi(t(X),t(X))X)
$$
therefore $f(t(X))=t(\Phi(t(X),t(X))X)=t\Phi(t(X),t(X))$. It follows that
$$f(x)=t\Phi(x,x)=t\phi_x(x).$$
% $ \sigma(\Phi(x,x)=f(x)$, for all $x\in M$, or $$\phi_x(x)=f(x).$$ 
%If $f_x=\beta\circ\phi_x$ then $f_x\in\mbox{Diff }(I)$ and $f_x(x)=f(x)$. 
The inverse is given by 
$$
\sigma^{-1}(y,Y)=(f^{-1}(y),\phi^{-1}_{f^{-1}(y)}(tY)Y).
$$
As a special case
\begin{equation*}\label{fifmenosum}
t(\phi_{f^{-1}(x)})^{-1}(x)=f^{-1}(x).
\end{equation*}

Let's denote by $\mathcal{J}$ the set of (local) right invariant diffeomorphisms of $I\times G$ that are $\rho_1$-projectable and preserve $B(G)$.  We have naturally the application
\begin{equation}\label{sigma}
\begin{array}{rcl}	\mathcal{J}&\rightarrow &\Gamma(G_k) \\
\sigma	&\mapsto  &\sigma_k
	\end{array},
\end{equation}
where $\sigma_k(x)=j^k_x\phi_x$, $x\in I$. If $\sigma'\in \mathcal{J}$, with $\sigma'=(f',\Phi')$, then
$$\begin{array}{rcl}	(\sigma'.\sigma)(x,y)&= &\sigma'(f(x),\phi_x(y)) \vspace{.2cm}\\
	&=  &(f'(f(x)),\phi'_{f(x)}(\phi_x(y))\vspace{.2cm}\\
	&=&((f'\circ f)(x),(\phi'_{f(x)}.\phi_x)(y)),
	\end{array}
$$
and from this it follows
$$(\sigma'\circ\sigma)_k(x)=j^k_x(\phi'_{f(x)}. \phi_x)=j^k_{f(x)}\phi'_{f(x)}.j^k_x\phi_x=\sigma'_k(f(x)).\sigma_k(x)=(\sigma'_k.\sigma_k)(x),$$
%$$(\sigma'\circ\sigma)_k(x)=j^k_x(\phi'_{f(x)}\circ \phi_x)=j^k_{f(x)}\phi'_{f(x)}.j^k_x\phi_x=\sigma'_k(f(x)).\sigma_k(x)=(\sigma'_k\circ\sigma_k)(x),$$
for each $x\in I$. So (\ref{sigma}) is a surjective morphism of groupoids. If $\phi$ is a rigth invariant diffeomorphism of $G$, let be $\tilde\phi\in \mathcal{J}$ given by
$$\tilde\phi(x,Y)=(t\phi(Y),\phi(tY)Y).$$
It is clear that
$$(\tilde\phi)_k=j^k\phi.$$

It follows from definitions of $\mathcal J$ and $\mathcal R$, that  is well defined the action
\begin{equation}\label{JR}
\begin{array}{rcl}
\mathcal{J}\times\mathcal{R}&\rightarrow &\mathcal R\\
(\sigma,v+\xi)&\mapsto &\sigma_*(v+\xi)
\end{array}.
\end{equation}
Then $\mathcal J$ acts on $\mathcal{V_R}$. Also, as $\tilde{\mathcal V}_{\mathcal R}$ is tangent to $B(G)$ and $B(G)$ is invariant by elements of $\mathcal J$, $\tilde{\mathcal V}_{\mathcal R}$ is invariant by $\mathcal J$.

%As $\sigma\in\mathcal{J}$ preserves $\mathcal{V_R}^{k+1}$, then the action of $\sigma$ on $\mathcal{R}$ quotients to an action on $\mathcal{R}/\mathcal{V_R}^{k+1}\cong\mathcal{T}\oplus J^k\mathcal{T}$.

\begin{po}\label{actionsigma}
Let be $\sigma\in\mathcal J$, $v\in\mathcal {H_R}$, $\xi\in\mathcal {V_R}$. We have:
\begin{enumerate}
\item [(i)] $(\sigma_*v)_k=f_*(v)+(j^1\sigma_k.v.\lambda^1\sigma_{k+1}^{-1}-\lambda^1\sigma_{k+1}.v.\lambda^1\sigma_{k+1}^{-1})$;
\item [(ii)]$(\sigma*\xi)_k=\lambda^1\sigma_{k+1}.\xi_k.\sigma_k^{-1}$;
\item [(iii)] $(\sigma*\tilde\xi)_k=f_*(\xi_H)+j^1\sigma_{k}.\xi_k.\sigma_k^{-1}$.
\end{enumerate}
\end{po}
\Pf If
$$\sigma(x,Y)=(f(x),\phi_x(tY)Y),$$
then
$$\sigma^{-1}(x,Y)=(f^{-1}(x),(\phi^{-1})_x(tY)Y),$$
where
$$(\phi^{-1})_x=(\phi_{f^{-1}(x)})^{-1}.$$

(i) Let be $v=\frac{d}{du}H_u|_{u=0}$ where $H_u(x,Y)=(h_u(x),Y)$. Therefore
$$
(\sigma_*v)=\frac{d}{du}(\sigma\circ H_u\circ \sigma^{-1})|_{u=0},
$$
or
$$
(\sigma_*v)(x,Y)=\frac{d}{du}\left((f\circ h_u\circ f^{-1})(x),(\phi_{(h_u\circ f^{-1})(x)}.(\phi_{f^{-1}(x)})^{-1})(tY)Y\right)|_{u=0}.
$$
Let be $\sigma\circ H_u\circ\sigma^{-1}=S_u\circ R_u$, where
$$
R_u(x,X)=(r_u(x),X)=\left((t\phi_{f^{-1}(x)}\circ h_u\circ f^{-1})(x),X\right)
$$
and
$$
S_u(y,X)=\left((f\circ h_u\circ f^{-1}\circ r_u^{-1})(y),(\phi_{(h_u\circ f^{-1}\circ r_u^{-1})(y)}.(\phi_{ (f^{-1}\circ r_u^{-1})(y)})^{-1})(tX)X\right).
$$
Therefore
$$
\begin{array}{rcl}
\frac{d}{du}S_u|_{u=0}(y,X)&=&\frac{d}{du}(f\circ h_u\circ f^{-1})|_{u=0}(y)-\frac{d}{du}r_u|_{u=0}(y)\\
&&+\frac{d}{du}(\phi_{(h_u\circ f^{-1}\circ r_u^{-1})(y)}. (\phi_{ (f^{-1}\circ r_u^{-1})(y)})^{-1})|_{u=0}(tX)X
\end{array}
$$
and
$$
\begin{array}{rcl}
\Upsilon_k(\frac{d}{du}S_u|_{u=0})(y)&=&(f_*v)(y)-\frac{d}{du}r_u|_{u=0}(y)\\
&&+j^k_y\left(\frac{d}{du}(\phi_{(h_u\circ f^{-1}\circ r_u^{-1})(y)}.(\phi_{ (f^{-1}\circ r_u^{-1})(y)})^{-1})\right)|_{u=0}.
\end{array}
$$
Let be $u\rightarrow x_u$ the family of trajectories   in $I$ defined by $x_u=r_u^{-1}(y)$. We can write
$$
\begin{array}{r}
j^k_y\left(\phi_{(h_u\circ f^{-1}\circ  r_u^{-1})(y)} .  (\phi_{ (f^{-1}\circ r_u^{-1})(y)})^{-1}\right)
=j^k_{(t\phi_{f^{-1}(x_u)}\circ h_u\circ f^{-1})(x_u)}\left(\phi_{(h_u\circ f^{-1})(x_u)}. (\phi_{ f^{-1}(x_u)})^{-1}\right)\\
=j^k_{(h_u\circ f^{-1})(x_u)}\phi_{(h_u\circ f^{-1})(x_u)}. j^k_{(t\phi_{f^{-1}(x_u)}\circ h_u\circ f^{-1})(x_u)}(\phi_{ f^{-1}(x_u)})^{-1}\\
=\sigma_k((h_u\circ f^{-1})(x_u)). j^k_{(t\phi_{f^{-1}(x_u)}\circ h_u\circ f^{-1})(x_u)}(\phi_{ f^{-1}(x_u)})^{-1}\\
\end{array}
$$
From this equality we get
$$
\begin{array}{rcl}
\Upsilon_k(\frac{d}{du}S_u|_{u=0})(y)&=&(f_*v)(y)-\frac{d}{du}r_u|_{u=0}(y)\\
&&+\frac{d}{du}\left[\sigma_k((h_u\circ f^{-1})(x_u)). j^k_{(t\phi_{f^{-1}(x_u)}\circ h_u\circ f^{-1})(x_u)}(\phi_{ f^{-1}(x_u)})^{-1}\right]_{u=0}\\
&=&(f_*v)(y)-\frac{d}{du}r_u|_{u=0}(y)\\
&&+\frac{d}{du}\left[\sigma_k((h_u\circ f^{-1})(y)).j^k_{(t\phi_{f^{-1}(y)}\circ h_u\circ f^{-1})(y)}(\phi_{ f^{-1}(y)})^{-1}\right]_{u=0}+\frac{d}{du}x_u|_{u=0}\\
&=&(f_*v)(y)-\frac{d}{du}r_u|_{u=0}(y)\\
&&+j^1_{f^{-1}(y)}\sigma_k. j^{k+1}_{y}(\phi_{ f^{-1}(y)})^{-1}.(\frac{d}{du}r_u|_{u=0}(y))+\frac{d}{du}x_u|_{u=0}\\
&=&(f_*v)(y)-\frac{d}{du}r_u|_{u=0}(y)\\
&&+j^1_{f^{-1}(y)}\sigma_k. (\sigma_{k+1})^{-1}(y).(\frac{d}{du}r_u|_{u=0}(y))+\frac{d}{du}x_u|_{u=0},
\end{array}
$$
since that $(\sigma_{k+1})^{-1}(y)=(\sigma_{k+1}(f^{-1}(y))^{-1}$. 

It follows from $r_u(x_u)=y$  that 
$$
\frac{d}{du}r_u|_{u=0}(y)+\frac{d}{du}x_u|_{u=0}=0.
$$
As $r(u)=(t\phi_{f^{-1}(x)}\circ h_u\circ f^{-1})(x)=
(t\phi_{f^{-1}(x)}\circ h_u\circ (t\phi_{f^{-1}(x)})^{-1})(x)=tj^{k}\phi_{f^{-1}(x)}\circ h_u\circ t(j^{k}\phi_{f^{-1}(x)})^{-1}(x)$, we obtain from proposition \ref{f41} that
%$$
%r(t)=tj^{k}\phi_{f^{-1}(x)}\circ h_t\circ (tj^{k}\phi_{f^{-1}(x)})^{-1}(x)
%$$
\begin{equation*}\label{2.12}
\frac{d}{du}r_u|_{u=0}=[\sigma_{k+1}.v.(\sigma_{k+1})^{-1}](x).
\end{equation*}
Then
$$
\begin{array}{rcl}
(\sigma_*v)_k(x)&=&\frac{d}{du}(S_u+H_u)|_{u=0}(x)\\
&=&(f_*v)(x)+j^1_{f^{-1}(x)}\sigma_k. (\sigma_{k+1})^{-1}(x).[\sigma_{k+1}.v.(\sigma_{k+1})^{-1}](x)-(\sigma_{k+1}.v.(\sigma_{k+1})^{-1})(x).
\end{array}
$$
We conclude 
$$
\begin{array}{rcl}
(\sigma_*v)_k=f_*v+j^1\sigma_{k}.v.(\sigma_{k+1})^{-1}-\sigma_{k+1}.v.(\sigma_{k+1})^{-1}
\end{array}
$$

(ii) Let be 
$$\xi=\frac{d}{du}V_u|_{u=0},$$
 where $V_u(x,Y)=(x,\eta^u_x(tY)Y)$, with $\eta^u_x\in\mathcal G$ for each $u$, and  $g_u(x)=t\eta^u_x(x)$. Then
$$(\sigma\circ V_u\circ \sigma^{-1})(x,Y)=(x,(\phi_{f^{-1}(x)}. \eta^u_{f^{-1}(x)}. (\phi_{f^{-1}(x)})^{-1})(tY)Y),$$
and 
$$(\sigma_*\xi)(x,Y)=\frac{d}{du}(\phi_{f^{-1}(x)}. \eta^u_{f^{-1}(x)}. (\phi_{f^{-1}(x)})^{-1})|_{u=0}(tY)Y
$$
Consequently
$$
\begin{array}{rcl}
\Upsilon_k(\sigma_*\xi)(x)&=&j^k_x(\frac{d}{du}(\phi_{f^{-1}(x)}. \eta^u_{f^{-1}(x)}.(\phi_{f^{-1}(x)})^{-1})|_{u=0})\vspace{.2cm}\\
&=&\frac{d}{du}(j^k_{(g_u\circ f^{-1})(x)}\phi_{f^{-1}(x)}. j^k_{f^{-1}(x)}\eta^u_{f^{-1}(x)}. j^k_x(\phi_{f^{-1}(x)})^{-1})|_{u=0}\vspace{.2cm}\\
&=&\frac{d}{du}\left((\tilde\phi_{f^{-1}(x)})_k((g_u\circ f^{-1})(x)). j^k_{f^{-1}(x)}\eta^u_{f^{-1}(x)}.((\tilde\phi_{f^{-1}(x)})^{-1})_k(x)\right)|_{u=0}\vspace{.2cm}\\
&=&\frac{d}{du}\left(j^1_{f^{-1}(x)}(\tilde\phi_{f^{-1}(x)})_k. j^k_{f^{-1}(x)}\eta^u_{f^{-1}(x)}.((\tilde\phi_{f^{-1}(x)})^{-1})_k(x)\right)|_{u=0}\vspace{.2cm}\\
&=&\lambda^1(\sigma_{k+1}(f^{-1}(x))).\xi_k(f^{-1}(x)).\sigma_k^{-1}(x),
\end{array}
$$
since that
$$j^1_{f^{-1}(x)}(\tilde\phi_{f^{-1}(x)})_k=j^1_{f^{-1}(x)}j^k\phi_{f^{-1}(x)}=\lambda^1(j^{k+1}_{f^{-1}(x)}\phi_{f^{-1}(x)})=\lambda^1(\sigma_{k+1}(f^{-1}(x))).$$
We proved
$$(\sigma*\xi)_k=\lambda^1\sigma_{k+1}.\xi_k.\sigma_k^{-1}.$$

(iii) Let be, as in (ii),  $\xi=\frac{d}{du}V_u|_{u=0}$, where $V_u(x,Y)=(x,\eta^u_x(tY)Y)$, with $\eta^u_x$  in  $\mathcal G$  for each $u$, and  $g_u(x)=t\eta^u_x(x)$. Then 
$$\tilde \xi=\frac{d}{du}\tilde V_u|_{u=0},$$
 where $\tilde V_u(x,Y)=(g_u(x),\eta^u_x(tY)Y)$.
Therefore
$$(\sigma\circ\tilde V_u\circ \sigma^{-1})(x,Y)=(f\circ g_u\circ f^{-1}(x),(\phi_{g_u\circ f^{-1}(x)}.\eta^u_{f^{-1}(x)}.(\phi_{f^{-1}(x)})^{-1})(tY)Y),$$
and 
$$(\sigma_*\tilde\xi)(x,Y)=f_*\xi_H(x)+\frac{d}{du}(\phi_{g_u\circ f^{-1}(x)}. \eta^u_{f^{-1}(x)}. (\phi_{f^{-1}(x)})^{-1})|_{u=0}(tY)Y
$$
Projecting, we obtain
$$
\begin{array}{rcl}
\Upsilon_k(\sigma_*\tilde\xi)(x)&=&f_*\xi_H(x)+j^k_x(\frac{d}{du}(\phi_{g_u\circ f^{-1}(x)}. \eta^u_{f^{-1}(x)}.(\phi_{f^{-1}(x)})^{-1})|_{u=0})\vspace{.15cm}\\
&=&f_*\xi_H(x)+\frac{d}{du}(j^k_{(g_u\circ f^{-1})(x)}\phi_{g_u\circ f^{-1}(x)}. j^k_{f^{-1}(x)}\eta^u_{f^{-1}(x)}. j^k_x(\phi_{f^{-1}(x)})^{-1})|_{u=0}\vspace{.15cm}\\
&=&f_*\xi_H(x)+\frac{d}{du}\left(\sigma_k((g_u\circ f^{-1})(x)). j^k_{f^{-1}(x)}\eta^u_{f^{-1}(x)}.(\sigma_k)^{-1}(x)\right)|_{u=0}\vspace{.15cm}\\
&=&f_*\xi_H(x)+j^1_{f^{-1}(x)}\sigma_{k}.\xi_k(f^{-1}(x)).\sigma_k^{-1}(x),
\end{array}
$$
so
$$(\sigma_*\tilde\xi)_k=f_*\xi_H+j^1\sigma_{k}.\xi_k.\sigma_k^{-1}.$$
Observe this formula depends only of $\sigma_k$.

We can give another proof combining (i) and (ii):
$$
\begin{array}{rcl}
(\sigma*\tilde\xi)_k&=&(\sigma*\xi_H)_k+(\sigma_*\xi)_k\\
&=&\left(f_*\xi_H+j^1\sigma_{k}.\xi_H.(\sigma_{k+1})^{-1}-\lambda^1\sigma_{k+1}.\xi_H.(\sigma_{k+1})^{-1}\right)+\left(\lambda^1\sigma_{k+1}.\xi_k.\sigma_k^{-1}\right)\\
&=&f_*\xi_H+j^1\sigma_{k}.\xi_H.(\sigma_{k+1})^{-1}+\sigma_{k+1}.(\xi_k-\xi_H).(\sigma_{k+1})^{-1}\\
\end{array}
$$
As $t_*(\xi_k-\xi_H)=0$, we get $j^1\sigma_k.(\xi_k-\xi_H)=\sigma_{k+1}.(\xi_k-\xi_H)$ so
$$\begin{array}{rcl}
(\sigma*\tilde\xi)_k&=&f_*\xi_H+j^1\sigma_{k}.\xi_H.(\sigma_{k+1})^{-1}+j^1\sigma_{k}.(\xi_k-\xi_H).(\sigma_{k+1})^{-1}\\
&=&f_*\xi_H+j^1\sigma_{k}.\xi_k.(\sigma_{k+1})^{-1}.
\end{array}
$$

\EPf

It follows from Proposition \ref{actionsigma}  that action (\ref{JR}) projects on an action $(\,\,)_*$:
\begin{equation}\label{GTG}
\begin{array}{rcl}
\mathcal G^{k+1}\times(\mathcal T\oplus \Gamma(J^k\mathfrak g))&\rightarrow &\mathcal T\oplus \Gamma(J^k\mathfrak g)\\
(\sigma_{k+1},v+\xi_k)&\mapsto & (\sigma_{k+1})_*(v+\xi_k)
\end{array}
\end{equation}
where
$$(\sigma_{k+1})_*(v+\xi_k)=f_*v+(j^1\sigma_k.v.\lambda^1\sigma_{k+1}^{-1}-\lambda^1\sigma_{k+1}.v.\lambda^1\sigma_{k+1}^{-1})+(\lambda^1\sigma_{k+1}.\xi_k.\sigma_k^{-1}).$$
This action verifies
\begin{equation}\label{colchactions}
\firstbra{(\sigma_{k+1})_*(v+\xi_k)}{(\sigma_{k+1})_*(w+\eta_k)}_k=(\sigma_k)_*(\firstbra{v+\xi_k}{w+\eta_k}_k).
\end{equation}
It follows from proposition \ref{actionsigma} and (\ref{GTG}) that  $(\sigma_{k+1})_*(\xi_k)(x)$ depends only on the value of $\sigma_{k+1}(x)$ at the point $x$ where $\xi$ is defined, and $(\sigma_{k+1})_*(v)(x)$ depends on the value of $\sigma_{k+1}$ on a curve tangent to $v(x)$. 

Item (iii) of proposition \ref{actionsigma} says that restriction to $J^k\tilde{\mathfrak g}$ of action (\ref {GTG}) is well defined:
\begin{equation*}\label{GG}
\begin{array}{rcl}
\mathcal G^{k}\times \Gamma(J^k\tilde{\mathfrak g})&\rightarrow & \Gamma(J^k\tilde{\mathfrak g})\\
(\sigma_{k},\tilde\xi_k)&\mapsto & (\sigma_{k})_*(\tilde\xi_k)
\end{array},
\end{equation*}
where 
\begin{equation*}\label{sigmatildexi}
(\sigma_k)_*(\tilde\xi_k)=f_*\xi_H+j^1\sigma_k.\xi_k.\sigma_k^{-1}.
\end{equation*}
In this case we get
$$\secondbra{(\sigma_k)_*\tilde\xi_k}{(\sigma_k)_*\tilde\eta_k}_k=(\sigma_k)_*\secondbra{\tilde\xi_k}{\tilde\eta_k}_k,
$$
and each $\nu_k(\sigma_k)_*\nu_k^{-1}$ acts as an automorphism of the Lie algebra sheaf $\Gamma(J^k{\mathfrak g})$:
$$\zerobra{\nu_k(\sigma_k)_*\nu_k^{-1}\xi_k}{\nu_k(\sigma_k)_*\nu_k^{-1}\eta_k}_k=\nu_k(\sigma_k)_*\nu_k^{-1}\zerobra{\xi_k}{\eta_k}_k
$$
%$$\zerobra{\nu(\sigma_k)_*\nu^{-1}\xi_k}{\nu(\sigma_k)_*\nu^{-1}\eta_k}_k=\nu^{-1}\nu(\sigma_k)_*\secondbra{\xi_k}{\eta_k}_k
%$$

%If $M$ and $M'$ are two manifolds of same dimension, we can define $$Q^k(M,M')=\{j^k_xf:f:U\subset M\rightarrow U'\subset M'\,\mbox{ is a diffeomorphism},x\in U\}.$$ 
%The groupoid $Q^k$ acts by the right on $Q^k(M,M')$, and $Q^{'k}=Q^k(M',M')$ acts on the left. Redoing the calculus of this subsection in this context, we obtain the analogous action of (\ref{GTG}): 
%\begin{equation}\label{GTGduplo}
%\begin{array}{rcl}
%\mathcal Q^{k+1}(M,M')\times(\mathcal T\oplus J^k\mathcal T)&\rightarrow &\mathcal T'\oplus J^k\mathcal T'\\
%(\sigma_{k+1},v+\xi_k)&\mapsto & (\sigma_{k+1})_*(v+\xi_k)
%\end{array},
%\end{equation}
%where $\mathcal Q^{k+1}(M,M')$ denotes the set of admissible sections of $ Q^{k+1}(M,M')$,
%\begin{equation}\label{actionsigmaformula}
%(\sigma_{k+1})_*(v+\xi_k)=f_*v+(j^1\sigma_k.v.\lambda^1\sigma_{k+1}^{-1}-\lambda^1\sigma_{k+1}.v.\lambda^1\sigma_{k+1}^{-1})+%(\lambda^1\sigma_{k+1}.\xi_k.\sigma_k^{-1}),
%\end{equation}
%and $f=\beta\sigma^{k+1}:M\rightarrow M'$.  This action also verifies (\ref{colchactions}).

%$\mathcal J^k\mathfrak g$

\subsection{The Lie algebra sheaf $\wedge( \check{\mathcal J}^\infty\mathfrak{g})^*\otimes ( \check{\cal J}^\infty\mathfrak g)$   }\label{las}

We denote by $\cal J^\infty\mathfrak g $ the projective limit of ${\cal J}^k\mathfrak g$, say,
$\cal J^\infty\mathfrak g$=$\lim\mbox{proj }{\cal J}^k\mathfrak g$,
and 
$$\check{\cal J}^\infty\mathfrak g=\mathcal T\oplus {\cal J}^\infty\mathfrak g.$$
As $T\oplus J^k\mathfrak g\cong TG^k|_I$, we have the identification of $\check{J}^\infty T$ with $\lim\mbox{proj } \Gamma(TG^k|_I)$, where $\Gamma(TG^k|_I)$ denotes the sheaf of germs  of local sections of the vector bundle $TG^k|_I\rightarrow I$.
From the fact that $\mathcal T\oplus {\cal J}^k\mathfrak g$ is a $\mathcal{O}_I$-module, we get  $\check{\cal J}^\infty\mathfrak g$ is a $\mathcal{O}_I$-module. In the following we use the notation
$$\check{\xi}=v+\lim\mbox{proj } \xi_k,\,\, \check{\eta}=w+\lim\mbox{proj } \eta_k\in\mathcal T\oplus {\cal J}^\infty\mathfrak g.$$
We define the first bracket in $\check{\cal J}^\infty\mathfrak g$ as:
\begin{equation*}\label{checkcolch}
\firstbra{\check{\xi}}{\check{\eta}}_\infty=\lim\mbox{proj }\firstbra{v+\xi_k}{w+\eta_k}_k
%(i(v)D\eta_{k+1}-i(w)D\xi_{k+1}+[\xi_{k+1},\eta_{k+1}]_{k+1}).
\end{equation*}
%The map 
%\begin{equation}\label{Rmap}
%\begin{array}{rcl}
%\mathcal R=\mathcal H\oplus\mathcal V&\rightarrow &\check{J}\mathcal A=\mathcal T\oplus J\mathcal A\\
%v+\xi&\mapsto&\check{\xi}=v+\lim\mbox{proj }\xi_k
%\end{array},
%\end{equation}
%is a surjective morphism of $\mathcal O_I$-modules, and from definition (\ref{checkcolch}) we obtain that the image of $[v+\xi,w+\eta]$ by the map (\ref{Rmap}) is $[\check{\xi},\check \eta]$.
With the bracket defined by (\ref{checkcolch}), $\check{\cal J}^\infty\mathfrak g$ is a 
\emph{Lie algebra sheaf}. Furthermore
\begin{equation*}\label{deriv}
\firstbra{\check\xi}{f\check \eta}_\infty=v(f)\check \eta+f\firstbra{\check\xi}{\check\eta}_\infty.
\end{equation*}
We extend now, as in \cite{Ma1}, \cite{Ma2} or \cite{KS}, the bracket on $\check{\cal J}^\infty\mathfrak g$ to a Nijenhuis bracket (see \cite{FN}) on $\wedge(\check{\cal J}^\infty\mathfrak g)^*\otimes (\check{\cal J}^\infty\mathfrak g)$, where 
\begin{equation*}\label{limind}
(\check{\cal J}^\infty\mathfrak g)^*=\lim\mbox{ind  }(\check{\cal J}^k\mathfrak g)^*.
\end{equation*}
We introduce the exterior differential $d$ on $\wedge(\check{\cal J}^\infty\mathfrak g)^*$, by:

(i) if $f\in\mathcal O_I$, then $df\in(\check{\cal J}^\infty\mathfrak g)^*$ is defined by
\begin{equation*}\label{df}
<df,\check\xi>=v(f).
\end{equation*}

(ii) if $\omega\in(\check{\cal J}^\infty\mathfrak g)^*$, then $d\omega\in\wedge^{2}(\check{\cal J}^\infty\mathfrak g)^*$ is defined by
\begin{equation*}\label{domega}
<d\omega,\check\xi\wedge\check\eta>=\mathfrak{L}(\check\xi)<\omega,\check\eta>-\mathfrak{L}(\check\eta)<\omega,\check\xi>-<\omega,\firstbra{\check\xi}{\check\eta}_\infty>,
\end{equation*}
where $\mathfrak{L}(\check\xi)f=<df,\check\xi>$.

We extend this operator to forms of any degree as a derivation of degree $+1$
$$d:\wedge^r(\check{\cal J}^\infty\mathfrak g)^*\rightarrow \wedge^{r+1}(\check{\cal J}^\infty\mathfrak g)^*.$$
The exterior differential $d$ is linear,
$$d(\omega\wedge\tau)=d\omega\wedge\tau+(-1)^{r}\omega\wedge d\tau,$$
for $\omega\in\wedge^r(\check{\cal J}^\infty\mathfrak g)^*$, and $d^2=0$.

Remember that  $(\rho_1)_*:\mathcal T\oplus {\cal J}^\infty\mathfrak g\rightarrow \mathcal T$  is the projection given by the decomposition in direct sum of $\check{\cal J}^\infty\mathfrak g=\mathcal{T}\oplus {\cal J}^\infty\mathfrak g$. (We could use, instead of $(\rho_1)_*$, the natural map $s_*:T\oplus J^k\mathfrak g\rightarrow T$, given by $s_*:TG_k|_I\rightarrow T$, and the identification (\ref{tqki})). Then $(\rho_1)^*:\mathcal T^*\rightarrow (\check{\cal J}^\infty\mathfrak g)^* $ and this map extends to $(\rho_1)^*:\wedge\mathcal T^*\rightarrow \wedge(\check{\cal J}^\infty\mathfrak g)^* $. If $\omega\in\wedge^r(\check{J}^\infty\mathcal T)^*$, then
$$<(\rho_1)^*\omega,\check\xi_1\wedge\cdots\wedge\check\xi_r>=<\omega,v_1\wedge\cdots
\wedge v_r>,$$
where $\check\xi_j=v_j+\xi_j,\,j=1,\cdots r$.
It follows that $d((\rho_1)^*\omega)=(\rho_1)^{*}(d\omega)$.
We identify $\wedge\mathcal T^*$ with its image in $\wedge(\check{\cal J}^\infty\mathfrak g)^*$ by $(\rho_1)^*$, and we write simply $\omega$ instead of $(\rho_1)^*\omega$. 

Let be $\textbf {u}=\omega\otimes\check\xi\in(\check{\cal J}^\infty\mathfrak g)^*\otimes (\check{\cal J}^\infty\mathfrak g)$, $\tau\in\wedge (\check{\cal J}^\infty\mathfrak g)^*$, with $\deg\omega=r$ and $\deg\tau=s$. We also define $\deg\textbf u =r$. Then we define the derivation of degree $(r-1)$
$$i(\textbf {u}):\wedge^s(\check{\cal J}^\infty\mathfrak g)^*\rightarrow \wedge^{s+r-1}(\check{\cal J}^\infty\mathfrak g)^*$$
by
\begin{equation}\label{i(u)}
i(\textbf {u})\tau=i(\omega\otimes\check\xi)\tau=\omega\wedge i(\check\xi)\tau
\end{equation}
and the Lie derivative
$$\mathfrak{L}(\textbf {u}):\wedge^r(\check{J}\mathfrak g)^*\rightarrow\wedge^{r+s}(\check{J}\mathfrak g)^*$$
by
\begin{equation}\label{theta(u)}
\mathfrak{L}(\textbf {u})\tau=i(\textbf {u})d\tau+(-1)^rd(i(\textbf {u})\tau),
\end{equation}
which is a derivation of degree $r$. If $\textbf {v}=\tau\otimes\check\eta$
we define 
\begin{equation}\label{buv}
[\textbf {u},\textbf {v}]=[\omega\otimes\check\xi,\tau\otimes\check\eta]=\omega\wedge\tau\otimes\firstbra{\check\xi}{\check\eta}_\infty+\mathfrak{L}(\omega\otimes\check\xi)\tau\otimes\check\eta-(-1)^{rs}\mathfrak{L}(\tau\otimes\check\eta)\omega\otimes\check\xi.
\end{equation}
A straightforward calculation shows that:
\begin{equation*}\label{bummed}
[\textbf {u},\tau\otimes\check\eta]=\mathfrak{L}(\textbf{u})\tau\otimes\check\eta+(-1)^{rs}\tau\wedge[\textbf u,\check\eta]-(-1)^{rs+s}d\tau\wedge i(\check\eta)\textbf{u},
\end{equation*}
where $i(\check \eta)\textbf u=i(\check\eta)(\omega\otimes\check\xi)=i(\check\eta)\omega\otimes\check\xi$. On verify that
$$
[\textbf {u},\textbf v]=-(-1)^{rs}[\textbf {v},\textbf u]
$$
and
\begin{equation}\label{brauv}
[\textbf {u},[\textbf v,\textbf w]]=[[\textbf u,\textbf v],\textbf w]+(-1)^{rs}[\textbf v,[\textbf u,\textbf w]],
\end{equation}
where $\deg\textbf u =r$, $\deg\textbf v =s$. 

With this bracket, $\wedge (\check{\cal J}^\infty\mathfrak g)^*\otimes (\check{\cal J}^\infty\mathfrak g)$ is a \emph{Lie algebra sheaf}.
Furthermore, if
\begin{equation*}\label{dtutv}
[\mathfrak{L}(\textbf u),\mathfrak{L}(\textbf v)]=\mathfrak{L}(\textbf u)\mathfrak{L}(\textbf v)-(-1)^{rs}\mathfrak{L}(\textbf v)\mathfrak{L}(\textbf u)
\end{equation*}
then
\begin{equation*}\label{tutv}
[\mathfrak{L}(\textbf u),\mathfrak{L}(\textbf v)]=\mathfrak{L}([\textbf u,\textbf v]).
\end{equation*}

In particular, we have the following formulas:
\begin{po}\label{p31}
If $\mbox{\bf u} \, ,\,\mbox{\bf v} \in (\check{\cal J}^\infty\mathfrak g)^*\otimes (\check{\cal J}^\infty\mathfrak g)$, $\omega\in(\check{\cal J}^\infty\mathfrak g)^*$, $\check\xi,\check\eta\in \check{\cal J}^\infty\mathfrak g$, then:
 
$$\begin{array}{lrcl}
(i)&<\mathfrak{L}(\mbox{\bf u})\omega,\check\xi\wedge\check\eta>&=&\mathfrak{L}(i(\check\xi)\mbox{\bf u})<\omega,\check\eta>-\mathfrak{L}(i(\check\eta)\mbox{\bf u})<\omega,\check\xi>\vspace{.2cm}\\ 
&&&-<\omega,\firstbra{i(\check\xi)\mbox{\bf u}}{\check\eta}_\infty+\firstbra{\check\xi}{i(\check\eta)\mbox{\bf u}}_\infty-i(\firstbra{\check\xi}{\check\eta}_\infty)\mbox{\bf u}>
\vspace{.3cm}\\
(ii)&i(\check\xi)[\mbox{\bf u},\check\eta]&=&\firstbra{i(\check\xi)\mbox{\bf u}}{\check\eta}_\infty-i(\firstbra{\check\xi}{\check\eta}_\infty)\mbox{\bf u},\vspace{.3cm}\\
(iii)&<[\mbox{\bf u},\mbox{\bf v}],\check\xi\wedge\check\eta>&=&\firstbra{i(\check\xi)\mbox{\bf u}}{i(\check\eta)\mbox{\bf v}}_\infty-\firstbra{i(\check\eta)\mbox{\bf u}}{i(\check\xi)\mbox{\bf v}}_\infty-\vspace{.2cm}\\
&&&-i(\firstbra{i(\check\xi)\mbox{\bf u}}{\check\eta}_\infty-\firstbra{i(\check\eta)\mbox{\bf u}}{\check\xi}_\infty-i(\firstbra{\check\xi}{\check\eta}_\infty)\mbox{\bf u})\mbox{\bf v}-\vspace{.2cm}\\
&&&-i(\firstbra{i(\check\xi)\mbox{\bf v}}{\check\eta}_\infty-\firstbra{i(\check\eta)\mbox{\bf v}}{\check\xi}_\infty-i(\firstbra{\check\xi}{\check\eta}_\infty)\mbox{\bf v})\mbox{\bf u}.
\end{array}$$
\end{po}
\Pf It is a straightforward calculus applying the definitions.
\EPf

If we define the groupoid
\begin{equation*}\label{mathcalG}
\mathcal G_{\infty}=\lim\mbox{proj }\mathcal G_k,
\end{equation*}
%From (\ref{sigma})  we get a morphism of groupoids:
%$$
%\begin{array}{rcl}
%\mathcal J&\rightarrow&\mathcal G_{\infty}\\
%\sigma&\mapsto&\lim\mbox{proj }\sigma_k,
%\end{array}
%
%$$
then for $\sigma=\lim\mbox{proj }\sigma_k\in \mathcal G_\infty$, we obtain, from (\ref{GTG}), the action
$$
\begin{array}{rcl}
\sigma_*:\check{\cal J}^\infty\mathfrak g&\rightarrow &\check{\cal J}^\infty\mathfrak g\\
\xi=v+\lim\mbox{proj }\xi_k&\mapsto&\sigma_*\xi=\lim\mbox{proj }(\sigma_{k+1})_*(v+\xi_k),
\end{array}
$$
so it is well defined
\begin{equation*}
\begin{array}{rcl}
\mathcal G_\infty\times\check{\cal J}^\infty\mathfrak g&\rightarrow&\check{\cal J}^\infty\mathfrak g\\
(\sigma,\xi)&\mapsto&\sigma_*\xi.
\end{array}.
\end{equation*}
It follows from (\ref{colchactions}) that $\sigma_*:\check{\cal J}^\infty\mathfrak g\rightarrow \check{\cal J}^\infty\mathfrak g$ is an \emph{automorphism of  Lie algebra sheaf}.
% $\check J\mathcal A$, as we can see from (\ref{colchactions}) and from definitions.

Given $\sigma\in\mathcal G_{\infty}$, $\sigma$ acts on $\wedge (\check{\cal J}^\infty\mathfrak g)^*$:
$$
\begin{array}{rcl}
\sigma^*:\wedge (\check J\mathcal T)^*&\rightarrow&\wedge (\check J\mathcal T)^*\\
\omega&\mapsto&\sigma^*\omega
\end{array},
$$
where, if $\omega$ is a $r$-form,
\begin{equation}\label{sigmaomega}
<\sigma^*\omega,\check\xi_1\wedge\cdots\wedge\check\xi_r>=<\omega, \sigma_*^{-1}(\check\xi_1)\wedge\cdots\wedge\sigma_*^{-1}(\check\xi_r)>.
\end{equation}
%and $\check\xi_j=v_j+\xi_j,\,j=1,\cdots r$.
Consequently, $\mathcal G_{\infty}$ acts on $\wedge(\check{\cal J}^\infty\mathfrak g)^*\otimes (\check{\cal J}^\infty\mathfrak g)$:
\begin{equation*}\label{GwJA}
\begin{array}{rcl}
\mathcal G_{\infty}\times(\wedge(\check{\cal J}^\infty\mathfrak g)^*\otimes (\check{\cal J}^\infty\mathfrak g))&\rightarrow&\wedge(\check{\cal J}^\infty\mathfrak g)^*\otimes (\check{\cal J}^\infty\mathfrak g)\\
(\sigma,\textbf u)&\mapsto&\sigma_*\textbf u
\end{array},
\end{equation*}
where
\begin{equation}\label{sigmastar}
\sigma_*\textbf u=\sigma_*(\omega\otimes\check\xi)=\sigma^*(\omega)\otimes\sigma_*(\check\xi).
\end{equation}
The action of $\sigma_*$ is an \emph{automorphism  of the Lie algebra sheaf} $\wedge(\check{\cal J}^\infty\mathfrak g)^*\otimes (\check{\cal J}^\infty\mathfrak g)$, i.e.,
$$[\sigma_*\textbf u,\sigma_*\textbf v]=\sigma_*[\textbf u,\textbf v].$$

\section{The first linear and non-linear Spencer complex}\label{fnsc}
In this section we will study the sub sheaf $\wedge\mathcal T^*\otimes \cal J^{\infty}\mathfrak g$ and introduce linear and non-linear Spencer complexes. Principal references are \cite{Ma1}, \cite{Ma2} and \cite{KS}.
\begin{po}\label{p41}
The sheaf $\wedge\mathcal T^*\otimes \cal J^{\infty}\mathfrak g$ is a Lie algebra  sub sheaf of $\wedge(\check {\cal J}^{\infty}\mathfrak g)^*\otimes (\check {\cal J}^{\infty}\mathfrak g)$, and
$$[\omega\otimes\xi,\tau\otimes \eta]=\omega\wedge\tau\otimes\firstbra{\xi}{\eta}_\infty,$$
where $\omega, \tau\in\wedge\mathcal T^*$, $\xi, \eta\in \cal J^{\infty}\mathfrak g$. 
\end{po}
\Pf Let be $\textbf u=\omega\otimes\xi\in\wedge\mathcal T^*\otimes \cal J^{\infty}\mathfrak g$. For any $\tau\in\wedge\mathcal T^*$,  $i(\xi)\tau=0$, then, applying (\ref{i(u)}) we obtain $i(\textbf u)\tau=0$, and by (\ref{theta(u)}), $\mathfrak{L}(\textbf u)\tau=0$. So (\ref{buv})
implies $[\omega\otimes\xi,\tau\otimes \eta]=\omega\wedge\tau\otimes\firstbra{\xi}{\eta}_\infty$.
\EPf

Let be the \emph{fundamental form}
$$\chi\in (\check {\cal J}^{\infty}\mathfrak g)^*\otimes (\check {\cal J}^{\infty}\mathfrak g)$$
defined by
$$i(\check\xi)\chi=(\rho_1)_*(\check\xi)=v,$$
where $\check\xi=v+\xi\in \mathcal T\oplus \cal J^{\infty}\mathfrak g$. In another words, $\chi$ is the projection of $\check {\cal J}^{\infty}\mathfrak g$ on $\mathcal T$, parallel to $\cal J^{\infty}\mathfrak g$.

If $\textbf u=\lim \textbf u_k$, we define $D\textbf u=\lim D\textbf u_k$.
\begin{po}\label{p42}
If $\omega\in\wedge\mathcal T^*$, and $\mbox{\bf u}\in\wedge\mathcal T^*\otimes \cal J^{\infty}\mathfrak g$, then:
\begin{enumerate}
\item [(i)] $\mathfrak{L}(\chi)\omega=d\omega$;
\item [(ii)] $[\chi,\chi]=0$;
\item [(iii)] $[\chi,\mbox{\bf u}]=D\mbox{\bf u}$.
\end{enumerate}
\end{po}
\Pf Let be $\check\xi=v+\xi,\,\,\check\eta=w+\eta\in\mathcal T\oplus \cal J^{\infty}\mathfrak g$.

(i) As $\mathfrak{L}(\chi)$ is a derivation of degree 1, it is enough to prove (i) for $0$-forms $f$ and $1$-forms $\omega\in(\check {\cal J}^{\infty}\mathfrak g)^*$. From (\ref{theta(u)}) we have $\mathfrak{L}(\chi)f=i(\chi)df=df$. It follows from proposition \ref{p31}(i) that
$$\begin{array}{rcl}
<\mathfrak{L}(\chi)\omega,\check\xi\wedge\check\eta>&=&\mathfrak{L}(v)<\omega,\check\eta>-\mathfrak{L}(w)<\omega,\check\xi>-<\omega,\firstbra{v}{\check\eta}_\infty+\firstbra{\check\xi}{w}_\infty-\chi(\firstbra{\check\xi}{\check\eta}_\infty)>\\
&=&\mathfrak{L}(v)<\omega, w>-\mathfrak{L}(w)<\omega, v>-<\omega,[v,w]+[v,w]-[v,w]>\\
&=&<d\omega,\check\xi\wedge\check\eta>.
\end{array}
$$

(ii) Applying proposition \ref{p31} (iii), we obtain
$$\begin{array}{rcl}
<\frac 1 2[\chi,\chi],\check\xi\wedge\check\eta>&=&[v,w]-i(\firstbra{v}{\check\eta}_\infty-\firstbra{w}{\check\xi}_\infty-\rho_1\firstbra{\check\xi}{\check\eta}_\infty)\chi\\
&=&[v,w]-([v,w]-[w,v]-[v,w])=0.
\end{array}
$$

(iii) It follows from  (\ref{buv}) that, for $\textbf u=\omega\otimes \xi$, 
$$
\begin{array}{rcl}
[\chi,\textbf u]&=&\mathfrak{L}(\chi)\omega\otimes\xi+(-1)^r\omega\wedge [\chi,\xi]-(-1)^{2r}d\omega\wedge i(\xi)\chi\\
&=&d\omega\otimes\xi+(-1)^r\omega\wedge[\chi,\xi].
\end{array}
$$
As $D$ is characterized by proposition \ref{propriedadesD}, it is enough to prove $[\chi,\xi]=D\xi$. It follows from propositions \ref{colchetek} (ii) and \ref{p31} (ii) that
$$i(\check\eta)[\chi,\xi]=\firstbra{i(\check\eta)\chi}{\xi}_\infty-i(\firstbra{\check\eta}{\xi}_\infty)\chi=\firstbra{w}{\xi}_\infty=i(\check\eta)D\xi.$$
\EPf

If $\textbf u,\textbf v\in \wedge\mathcal T^*\otimes \cal J^{\infty}\mathfrak g$, with $\deg \textbf u=r$, $\deg \textbf v=s$, then we get from (\ref{brauv}) and proposition \ref{p42} (iii) that
$$D[\textbf u,\textbf v]=[D\textbf u,\textbf v]+(-1)^r[\textbf u,D\textbf v],$$
and
$$[\chi,[\chi,\textbf u]]=[[\chi,\chi],\textbf u]-[\chi,[\chi,\textbf u]]=-[\chi,[\chi,\textbf u]],$$
therefore, $D^2\textbf u=0$, or
$$D^2=0.$$
Then it is well defined \emph{the first linear Spencer complex},
$$0\rightarrow \Gamma(\mathfrak g)\stackrel{j^\infty}{\rightarrow}{\cal J}^{\infty}\mathfrak g\stackrel{D}{\rightarrow}{\mathcal T}^*\otimes \cal J^{\infty}{\mathfrak g}\stackrel{D}{\rightarrow}$$
$${\wedge}^2{\mathcal T}^*\otimes {\cal J}^{\infty}{\mathfrak g}\stackrel{D}
{\rightarrow}\cdots
\stackrel{D}{\rightarrow}{\wedge}^m{\mathcal T}^*\otimes {\cal J}^{\infty}\mathfrak g\rightarrow 0,$$
where $\dim T=m$. This complex projects on
$$0\rightarrow \Gamma(\mathfrak g)\stackrel{j^k}{\rightarrow} {\cal J}^{k}\mathfrak g\stackrel{D}{\rightarrow}{\mathcal T}^*\otimes {\cal J}^{k-1}\mathfrak g\stackrel{D}{\rightarrow}{\wedge^2}{\mathcal T}^*\otimes {\cal J}^{k-2}\mathfrak g\stackrel{D}{\rightarrow}\cdots 
 \stackrel{D}{\rightarrow}\wedge^{m}{\mathcal T}^*\otimes {\cal J}^{k-m}\mathfrak g\rightarrow 0.$$
This complex is exact (see \cite {Ma1}, \cite{Ma2}, \cite{KS}).

Let be $\gamma^k$ the kernel of $\pi_k:J^k\mathfrak g\rightarrow J^{k-1}\mathfrak g$. Denote by $\delta$ the restriction of $D$ to $\gamma^k$.   It follows from proposition \ref{propriedadesD}(ii) that $\delta$ is $\mathcal{O}_I$-linear and $\delta:\gamma^k\rightarrow T^*\otimes \gamma^{k-1}$. 
This map is injective, in fact, if  $\xi\in\gamma^k$, then by (\ref{kernelpi}), $\delta\xi=-\lambda^1(\xi)$ is injective. 
As
$$i(v)D(i(w)D\xi)-i(w)D(i(v)D\xi)-i([v,w])D\pi_{k-1}\xi=0,$$
for $v,w\in\mathcal T$, $\xi\in \gamma^k\subset {\cal J}^k\mathfrak g$, we obtain that $\delta$ is symmetric, $i(v)\delta(i(w)\delta\xi)=i(w)\delta(i(v)\delta\xi).$
Observe that we get the map
$$\iota:\gamma^k\rightarrow S^2T^*\otimes\gamma^{k-2}$$
defined by
$i(v,w)\iota(\xi)=i(w)\delta(i(v)\delta\xi),$ and if we go on, we obtain the isomorphism
$$\gamma^k\cong S^kT^*\otimes J^0\mathfrak g,$$
where, given   basis 
$e_1,\cdots,e_m\in T$ , $f_1,\cdots,f_r\in \mathfrak g$
with  the dual basis $e^1,\cdots,e^m\in T^*$,
we obtain the basis
\begin{equation*}\label{f}
f^{k_1,k_2,\cdots,k_m}_l=\frac{1}{k_1!k_2!\cdots k_m!}(e^1)^{k_1}(e^2)^{k_2}\cdots (e^m)^{k_m}\otimes j^0f_l
\end{equation*}  
of $S^kT^*\otimes J^0\mathfrak g$, where $k_1+k_2+\cdots+k_m=k$,  $k_1,\cdots,k_m\geq 0$ and $l=1,\cdots,r$. In this basis 
$$\delta(f^{k_1,k_2,\cdots,k_m}_l)=-\Sigma_{i=1}^me^i\otimes f^{k_1,\cdots,k_{i-1},k_i-1,k_{i+1},\cdots,k_m}_l.$$
From the linear Spencer complex, we obtain the exact sequence of morphisms of vector bundles
\begin{equation*}\label{delta}
0\rightarrow \gamma^k\stackrel{\delta}{\rightarrow} T^*\otimes \gamma^{k-1}\stackrel{\delta}{\rightarrow}\wedge^2\mathcal T^*\otimes \gamma^{k-2} \stackrel{\delta}{\rightarrow}\cdots\stackrel{\delta}{\rightarrow}\wedge^m T^*\otimes \gamma^{k-m} \rightarrow 0.
\end{equation*}

Let's now introduce the \emph{first nonlinear Spencer operator} $\mathcal D$. The ``finite'' form $\mathcal D$ of the linear Spencer  operator $D$ is defined by
\begin{equation*}\label{dxi}
\mathcal D\sigma=\chi-\sigma_*^{-1}(\chi),
\end{equation*}
where $\sigma\in\mathcal G_{\infty}$.
\begin{po}\label{p43}
The operator $\mathcal D$ take values in $\mathcal T^*\otimes {\cal J}^\infty\mathfrak g$, so
$$\mathcal D:\mathcal G_{\infty}\rightarrow \mathcal T^*\otimes {\cal J}^\infty\mathfrak g,$$
and
\begin{equation}\label{ivds}
i(v)(\mathcal D\sigma)_k=\lambda^1\sigma_{k+1}^{-1}.j^1\sigma_k.v-v,
\end{equation}
where $\sigma  = \lim \emph{proj }\sigma_k\in\mathcal G_{\infty}$.
\end{po}
\Pf  Applying (\ref{sigmaomega}) and (\ref{sigmastar}), it follows for  $\xi\in {\cal J}^\infty\mathfrak g$,
$$i(\xi)\mathcal D\sigma=i(\xi)\chi-\sigma^{-1}_*(i(\sigma_*(\xi))\chi)=0,$$
and for $v\in\mathcal T$,
$$i(v)\mathcal D\sigma=i(v)\chi-\sigma^{-1}_*(i(\sigma_*(v))\chi)=v-\sigma_*^{-1}(f_*v),$$
where $f=t\circ\sigma$. By proposition \ref{actionsigma} (i), 
$$i(v)(\mathcal D\sigma)_k=v-(f_*^{-1}(f_*v)+j^1\sigma_k^{-1}.f_*v.\lambda^1\sigma_{k+1}-\lambda^1\sigma_{k+1}^{-1}.f_*v.\lambda^1\sigma_{k+1}).$$
Posing $v=\frac{d}{du}x_u|_{u=0}$, we obtain
\begin{equation}\label{por}
j^1\sigma_k.v.j^1\sigma_k^{-1}=\frac{d}{du}(\sigma_k(x_u).\sigma_k^{-1}(f(x_u)))|_{t=0}=\frac{d}{du}f(x_u)|_{u=0}=f_*v,
\end{equation}
and replacing this above, we get
$$
\begin{array}{rcl}
i(v)(\mathcal D\sigma)_k&=&-j^1\sigma_k^{-1}.(j^1\sigma_k.v.j^1\sigma_k^{-1}).\lambda^1\sigma_{k+1}+\lambda^1\sigma_{k+1}^{-1}.(j^1\sigma_k.v.j^1\sigma_k^{-1}).\lambda^1\sigma_{k+1}\vspace{.2cm}\\
&=&(\lambda^1\sigma_{k+1}^{-1}.j^1\sigma_k.v-v).j^1\sigma_k^{-1}.\lambda^1\sigma_{k+1}\vspace{.2cm}\\
&=&(\lambda^1\sigma_{k+1}^{-1}.j^1\sigma_k.v-v).\sigma_k^{-1}.\sigma_{k}\vspace{.2cm}\\
&=&\lambda^1\sigma_{k+1}^{-1}.j^1\sigma_k.v-v,
\end{array}
$$
since that $\lambda^1\sigma_{k+1}^{-1}.j^1\sigma_k.v-v$ is  $s$-vertical (cf. (\ref{actionvert})).
\EPf
\begin{co}
We have ${\mathcal D}\sigma=0$ if and only if $\sigma=j^{\infty}(\pi_0\circ\sigma)$, where $\pi_0:\mathcal G_{\infty} \rightarrow \mathcal  G$.
\end{co}
\begin{co}\label{sigmaD}
If $\sigma_{k+1}\in\mathcal{G}_{k+1}$, then
$$(\sigma_{k+1})_*(v)=f_*v+(\sigma_{k+1})_*(i(v)\mathcal{D}\sigma_{k+1}),$$
for $v\in \mathcal{T}$.
\end{co}
\Pf It follows from (\ref{GTG}) and proposition \ref{p43} that  
$$
\begin{array}{rcl}(\sigma_{k+1})_*(i(v)\mathcal{D}\sigma_{k+1})&=&\lambda^1\sigma_{k+1}.(i(v)\mathcal{D}\sigma_{k+1}).\lambda^1\sigma_{k+1}^{-1}=\lambda^1\sigma_{k+1}.(\lambda^1\sigma_{k+1}^{-1}.j^1\sigma_k.v-v).\lambda^1\sigma_{k+1}^{-1}\vspace{.2cm}\\
&=&j^1\sigma_k.v.\lambda^1\sigma_{k+1}^{-1}-\lambda^1\sigma_{k+1}.v.\lambda^1\sigma_{k+1}^{-1}=(\sigma_{k+1})_*(v)-f_*v.
\end{array}
$$
\EPf
\begin{po}\label{p44}
The operator $\mathcal D$ has the following properties:
\begin{enumerate}
\item [(i)] If $\sigma,\sigma'\in \mathcal G_{\infty}$,
$$\mathcal D(\sigma'\circ\sigma)=\mathcal D\sigma+\sigma_*^{-1}(\mathcal D\sigma').$$
In particular
$$\mathcal D \sigma^{-1}=-\sigma_*(\mathcal D\sigma).$$
\item [(ii)] If $\sigma\in\mathcal G_{\infty}$, $\mbox{\bf u}\in\wedge\mathcal T^*\otimes {\cal J}^\infty\mathfrak g$,
$$ \mathcal D(\sigma^{-1}_*\mbox{\bf u})=\sigma_*^{-1}(D\mbox{\bf u})+[\mathcal D\sigma,\sigma^{-1}_*\mbox{\bf u}].$$
\item [(iii)] If $\xi=\frac{d}{du}\sigma_u|_{u=0}$, with $\xi\in {\cal J}^\infty\mathfrak g$, and $\sigma_u\in\mathcal G_{\infty}$ is the 1-parameter group associated to $\xi$, then
$$D\xi=\frac{d}{du}\mathcal D\sigma_u|_{u=0}.$$
\end{enumerate}
\end{po}
\Pf 

(i) $$\mathcal D(\sigma'\circ\sigma)=\chi-\sigma_*^{-1}(\chi)+\sigma_*^{-1}(\chi-(\sigma')_*^{-1}(\chi))=\mathcal D\sigma+\sigma_*^{-1}(\mathcal D\sigma').$$

(ii) $$\begin{array}{rcl}
 D(\sigma_*^{-1}\textbf u)&=&\sigma_*^{-1}[\sigma_*\chi,\textbf u]=\sigma_*^{-1}[\chi-\mathcal D\sigma^{-1},\textbf u]\vspace{.3cm}\\
&=&\sigma_*^{-1}(D\textbf u)-[\sigma_*^{-1}(\mathcal D\sigma^{-1}),\sigma_*^{-1}\textbf u]=\sigma_*^{-1}(D\textbf u)+[\mathcal D\sigma,\sigma_*^{-1}\textbf u].
\end{array}$$

(iii) $$\frac{d}{du}\mathcal D\sigma_u|_{u=0}=-\frac{d}{du}(\sigma_u^{-1})_*(\chi)=-[\xi,\chi]=D\xi.$$
\EPf

Proposition \ref{p43} says that $\mathcal D$ is projectable:
\begin{equation*}\label{f42}
\begin{array}{rcl}
\mathcal D:\mathcal G_{k+1}&\rightarrow &\mathcal T^*\otimes {\cal J}^k\mathfrak g\\
\sigma_{k+1}&\mapsto &\mathcal D\sigma_{k+1}
\end{array},
\end{equation*}
where $$i(v)\mathcal D\sigma_{k+1}=\lambda^1\sigma_{k+1}^{-1}.j^1\sigma_k.v-v.$$ 
%(see formula (6.8) of \cite{Ma2})

It follows from $[\chi,\chi]=0$ that 
$$0=\sigma_*^{-1}([\chi,\chi])=[\sigma_*^{-1}(\chi),\sigma_*^{-1}(\chi)]=[\chi-\mathcal D\sigma,\chi-\mathcal D\sigma]=[\mathcal D\sigma,\mathcal D\sigma]-2D(\mathcal D\sigma),$$
therefore
\begin{equation}\label {f43}
D(\mathcal D\sigma)-\frac{1}{2}[\mathcal D\sigma,\mathcal D\sigma]=0.
\end{equation}
If we define the non linear operator 
\begin{equation*}\label{f44}
\begin{array}{rcl}
\mathcal D_1:\mathcal T^*\otimes J^\infty\mathfrak g&\rightarrow &\wedge^2\mathcal T^*\otimes J^\infty\mathfrak g\\
\textbf u&\mapsto & D\textbf u-\frac{1}{2}[\textbf u,\textbf u]
\end{array},
\end{equation*}
then we can write  (\ref{f43}) as
$$\mathcal D_1\mathcal D=0.$$
The operator $\mathcal D_1$ projects in order $k$ to
$$\mathcal D_1:\mathcal T^*\otimes {\mathcal J}^k\mathfrak g\rightarrow \wedge^2\mathcal T^*\otimes {\mathcal J}^{k-1}\mathfrak g,$$
where
$$\mathcal D_1\textbf u=D\textbf u-\frac{1}{2}[\textbf u,\textbf u]_k.$$
We define the \emph{first non-linear Spencer complex} by
$$1\rightarrow {\mathcal G}\stackrel{j^{k+1}}{\rightarrow}\mathcal G_{k+1}\stackrel{\mathcal D}{\rightarrow} \mathcal T^*\otimes {\mathcal J}^k\mathfrak g\stackrel{\mathcal D_1}{\rightarrow}\wedge^2\mathcal T^*\otimes {\mathcal J}^{k-1}\mathfrak g,$$
which is exact in $\mathcal G_{k+1}$.

Let be $p^1({\mathcal D}):J^1G_{k+1}\rightarrow  T^*\otimes J^k\mathfrak g$ the morphism associated to the differential operator ${\mathcal D}$. It follows from proposition \ref{p43} that 
$$
i(v)p^1({\mathcal D})(j^1_x\sigma_{k+1})=\lambda^1\sigma_{k+1}^{-1}(x).j^1_x\sigma_k.v-v
$$
where $\sigma_{k+1}$ is an admissible section of $G_{k+1}$ and $\sigma_k=\pi_k\sigma_{k+1}$.

\begin{po} The image of $G_{k,1}$ by $p^1({\mathcal D})$ is the set
$$
 B^{k,1}=\{X\in { T}^*\otimes  J^k\mathfrak g: v\in T\rightarrow t_*(X(v)+v)\in T \mbox{ is inversible}\}.
$$
\end{po}
\Pf It follows from (\ref{ivds}) and (\ref{por}) that 
$t_*(i(v){\mathcal D}\sigma_{k+1})=
\sigma_{k+1}^{-1}.f_*v.\sigma_{k+1}-v$, so $t_*(i(v){\mathcal D}\sigma_{k+1}+v)=\sigma_{k+1}^{-1}.f_*v.\sigma_{k+1}$ therefore ${\mathcal D}\sigma_{k+1}\in  B^{k,1}$. 

Conversely, let be $Y\in  B^{k,1}_x$ and consider $\bar Y\in  B^{k+1,1}_x$ such that $\mbox{id}^*\otimes\pi_k(\bar Y)=Y$. Let be the map $\Sigma:T_x\rightarrow T_{I_{k+1}(x)}G_k$ defined by $\Sigma(v)=i(v)\bar Y+v$.   There exists a local section 
$\sigma\in \mathcal{G}_{k+1}$ such that $\sigma(x)=I_{k+1}(x)$ and $\sigma_*(v)=\Sigma(v)$. Then 
$$
\begin{array}{rcl}
i(v)p^1({\mathcal D})(j^1_x\sigma)&=&\sigma^{-1}(x).j^1_x(\pi_k\sigma).v-v\\
&=&j^1_x(\pi_k\sigma).v-v\\
&=&\pi_ki(v)\Sigma-v\\
&=&(i(v) Y+v)-v\\
&=&i(v)Y.
\end{array}
$$
\EPf

%It is possible  to define the first nonlinear Spencer complex $\mathcal D$ for admissible sections of $Q^{\infty}(M,M')$ by:
%\begin{equation}\label{dxig}
%\mathcal D\sigma=\chi-\sigma_*^{-1}(\chi'),
%\end{equation}
%where $\sigma\in\mathcal G_{\infty}(M,M')$ and $\chi'\in (\check J^\infty\mathcal T')^*\otimes (\check J^\infty\mathcal T')$ is the fundamental form.
%\begin{po}%\label{p43}
%The operator $\mathcal D$ take values in $\mathcal T^*\otimes {\cal J}^\infty\mathfrak g$, so
%$$\mathcal D:\mathcal G_{\infty}(M,M')\rightarrow \mathcal T^*\otimes {\mathcal J}^\infty\mathfrak g,$$
%and the same formula of proposition \ref{p43} holds:
%$$i(v)(\mathcal D\sigma)_k=\lambda^1\sigma_{k+1}^{-1}.j^1\sigma_k.v-v,$$
%where $\sigma  = \lim \mbox{proj }\sigma_k\in\mathcal Q^{\infty}(M,M')$. Other properties can easily be generalised.

%%%%%%%%%%%%%%%%%%%%%%%
%%%%%%%%%%%%%%%%%%%%
%%%%%%%%%%%%%%%%
%%%%%%%%%%%%
%%%%%%%
\section{The second linear and non-linear Spencer complex}\label{snsc}

Consider the projection $\rho:\check J^\infty\mathfrak g\rightarrow T$ defined by $\rho(v+\xi)=v-t_*(\xi)$. The kernel of $\rho$ is $\tilde J\mathfrak g$. We can include $\wedge{\mathcal T}^*$ in $\check J^\infty\mathfrak g$ by the pullback for $\rho^*$. Therefore we denote by   $\wedge\tilde{\mathcal T^*}=\rho^*(\wedge{\mathcal T}^*)$.

\begin{po}\label{p42ii}
The sheaf $\wedge\tilde{\mathcal T}^*\otimes \tilde{\cal J}^{\infty}\mathfrak g$ is a Lie algebra  sub sheaf of $\wedge(\check {\cal J}^{\infty}\mathfrak g)^*\otimes (\check {\cal J}^{\infty}\mathfrak g)$.
and
%$$[\omega\otimes\xi,\tau\otimes \eta]=\omega\wedge\tau\otimes\firstbra{\xi}{\eta}_\infty,$$
%where $\omega, \tau\in\wedge\mathcal T^*$, $\xi, \eta\in \cal J^{\infty}\mathfrak g$. 
$$[\tilde\omega\otimes\tilde\xi,\tilde\tau\otimes \tilde\eta]=\tilde\omega\wedge\tilde\tau\otimes\secondbra{\tilde\xi}{\tilde\eta}_\infty,$$
where $\tilde\omega, \tilde\tau\in\wedge\tilde{\mathcal T}^*$, $\tilde\xi, \tilde\eta\in \tilde{\mathcal J}^{\infty}\mathfrak g$. 
\end{po}
\Pf Let be $\textbf u=\tilde\omega\otimes\tilde\xi\in\wedge\tilde{\mathcal T}^*\otimes \tilde{\mathcal J}^{\infty}\mathfrak g$. For any $\tilde\tau\in\wedge\tilde{\mathcal T}^*$,  $i(\tilde\xi)\tilde\tau=0$, then, applying (\ref{i(u)}) we obtain $i(\textbf u)\tilde\tau=0$, and by (\ref{theta(u)}), $\mathfrak{L}(\textbf u)\tilde\tau=0$. So (\ref{buv})
implies $[\tilde\omega\otimes\tilde\xi,\tilde\tau\otimes \tilde\eta]=\tilde\omega\wedge\tilde\tau\otimes\secondbra{\xi}{\eta}_\infty$.
\EPf

%Let be $\overline \chi\in (\check {\cal J}^{\infty}\mathfrak g)^*\otimes (\check {\mathcal J}^{\infty}\mathfrak g)$ the projection of $\check {\cal J}^{\infty}\mathfrak g$ on $\mathcal{T}$ parallelly to  $\tilde {\mathcal J}^{\infty}\mathfrak g$, i.e.,
%$$
%\overline\chi(\check\xi)=\overline\chi(v+\xi)=v-t_*(\xi).
%$$ 
%If $\omega\in\wedge^r{\mathcal T}^*$  we define $\tilde\omega\in\wedge^r(\check {\mathcal J}^{\infty}\mathfrak g)^*$   by
%$$
%\langle \tilde\omega,\check\xi_1\wedge\cdots\wedge\check\xi_r\rangle=\langle\omega,i(\check\xi_1)\overline\chi\wedge\cdots\wedge i(\check\xi_r)\overline\chi\rangle.
%$$
%Let be 
%$$\wedge\tilde{\mathcal T}^*=\{\tilde\omega\in \wedge(\check {\mathcal J}^{\infty}\mathfrak g)^*| \omega\in\wedge\mathcal T^* \}$$
 %and we introduce
 %$$\begin{array}{rcl}
 %\overline d:\wedge\tilde{\mathcal T}^*&\rightarrow &\wedge\tilde{\mathcal T}^*\\
 %\tilde\omega&\rightarrow&\overline d\tilde\omega=\widetilde{d\omega}
 %\end{array}$$

\begin{co}
$\wedge\tilde{\mathcal T}^*\otimes \tilde{\mathcal J}^{k}\mathfrak g$ is a sheaf in Lie algebras for $k\geq 1$.
\end{co}

\Pf The bracket $\secondbra{}{}_\infty$ projects to $\secondbra{}{}_k$ in  $\tilde{\mathcal J}^{k}\mathfrak g$. Therefore $\wedge\tilde{\mathcal T}^*\otimes \tilde{\mathcal J}^{k}\mathfrak g$ is well defined as a sheaf.
\EPf

Let be the \emph{fundamental form}
$$\overline\chi\in (\check {\cal J}^{\infty}\mathfrak g)^*\otimes (\check {\cal J}^{\infty}\mathfrak g)$$
defined by
$$i(\check\xi)\overline\chi=(\rho_1-t)_*(\check\xi)=v-t_*(\xi),$$
where $\check\xi=v+\xi\in \mathcal T\oplus \cal J^{\infty}\mathfrak g$. In another words, $\overline\chi$ is the projection of $\check {\cal J}^{\infty}\mathfrak g$ on $\mathcal T$, parallel to $\tilde{\cal J}^{\infty}\mathfrak g$.

We define \emph{the second linear Spencer operator} $\overline D$ by
$$
\begin{array}{rcl}
\overline D:\wedge\tilde{\mathcal T}^*\otimes \tilde{\mathcal J}^{\infty}\mathfrak g&\rightarrow&\wedge\tilde{\mathcal T}^*\otimes \tilde{\mathcal J}^{\infty}\mathfrak g\\
\textbf u&\rightarrow &\overline D\textbf u=\nu ^{-1}D\nu\textbf u,
\end{array}
$$
where 
$$
\begin{array}{rcl}
\nu:\wedge\tilde{\mathcal T}^*\otimes \tilde{\mathcal J}^{\infty}\mathfrak g&\rightarrow&\wedge{\mathcal T}^*\otimes {\mathcal J}^{\infty}\mathfrak g\\
\tilde\omega\otimes \tilde\xi&\rightarrow &\omega\otimes\xi.
\end{array}
$$
We project this isomorphism $\nu$ in order $k$ to
\begin{equation*}%\label{nu}
\begin{array}{rcl}
\nu_k:\wedge\tilde{\mathcal T}^*\otimes \tilde{\mathcal J}^{k}\mathfrak g&\rightarrow&\wedge{\mathcal T}^*\otimes {\mathcal J}^{k}\mathfrak g\\
\tilde\omega\otimes \tilde\xi_k&\rightarrow &\omega\otimes\xi_k.
\end{array}
\end{equation*}
%If $\textbf u=\lim \textbf u_k$, we define $\overline D\textbf u=\lim D\textbf u_k$.
\begin{po}\label{p42iii}
If $\tilde\omega\in\wedge\tilde{\mathcal T}^*$, and $\mbox{\bf u}\in\wedge\tilde{\mathcal T}^*\otimes \tilde{\mathcal J}^{\infty}\mathfrak g$, then:
\begin{enumerate}
\item [(i)] $\mathfrak{L}(\overline\chi)\tilde\omega={d\tilde\omega}$;
\item [(ii)] $[\overline\chi,\overline\chi]=0$;
\item [(iii)] $[\overline\chi,\mbox{\bf u}]=\overline D\mbox{\bf u}$.
\end{enumerate}
\end{po}
\Pf Let be $\check\xi=v+\xi,\,\,\check\eta=w+\eta\in\mathcal T\oplus \cal J^{\infty}\mathfrak g$. 

(i) As $\mathfrak{L}(\overline\chi)$ is a derivation of degree 1, it is enough to prove (i) for $0$-forms $f$ and $1$-forms $\tilde\omega\in(\check {\mathcal J}^{\infty}\mathfrak g)^*$. From (\ref{theta(u)}) we have $\mathfrak{L}(\overline\chi)f=i(\overline\chi)df=\widetilde{df}$. It follows from proposition \ref{p31}(i) that
$$\begin{array}{rcl}
<\mathfrak{L}(\overline\chi)\tilde\omega,\check\xi\wedge\check\eta>&=&\mathfrak{L}(v-\xi_H)<\tilde\omega,\check\eta>-\mathfrak{L}(w-\eta_H)<\tilde\omega,\check\xi>\\
&&-<\omega,\firstbra{v-\xi_H}{\check\eta}_\infty+\firstbra{\check\xi}{w-\eta_H}_\infty-\overline\chi(\firstbra{\check\xi}{\check\eta}_\infty)>\\
&=&\mathfrak{L}(v-\xi_H)<\omega, w-\eta_H>-\mathfrak{L}(w-\eta_H)<\omega, v-\xi_H>\\
&&-<\omega,[v-\xi_H,w]-t_*(i(v-\xi_H)D\eta)+[v,w-\eta_H]+t_*(i(w-\eta_H)D\xi)>\\
&&+<\omega, [v,w]-t_*(i(v)D\eta-i(w)D\xi+\firstbra{\xi}{\eta}_\infty)>\\
&=&\mathfrak{L}(v-\xi_H)<\omega, w-\eta_H>-\mathfrak{L}(w-\eta_H)<\omega, v-\xi_H>\\
&&-<\omega,[v,w]-[\xi_H,w]-[v,\eta_H]+t_*(i(\xi_H)D\eta-i(\eta_H)D\xi+\firstbra{\xi}{\eta}_\infty)>\\
&=&\mathfrak{L}(v-\xi_H)<\omega, w-\eta_H>-\mathfrak{L}(w-\eta_H)<\omega, v-\xi_H>\\
&&-<\omega,[v,w]-[\xi_H,w]-[v,\eta_H]+t_*([\tilde\xi,\tilde\eta]_\infty>\\
&=&\mathfrak{L}(v-\xi_H)<\omega, w-\eta_H>-\mathfrak{L}(w-\eta_H)<\omega, v-\xi_H>\\
&&-<\omega,[v,w]-[\xi_H,w]-[v,\eta_H]+[\xi_H,\eta_H]>\\
&=&\mathfrak{L}(v-\xi_H)<\omega, w-\eta_H>-\mathfrak{L}(w-\eta_H)<\omega, v-\xi_H>\\
&&-<\omega,[v-\xi_H,w-\eta_H]>\\
&=&<d\omega,\rho_*(\check\xi)\wedge \rho_*(\check\eta)>=<d\tilde\omega,\check\xi\wedge\check\eta>.
\end{array}
$$

(ii) Applying proposition \ref{p31} (iii), we obtain
$$\begin{array}{rcl}
<\frac 1 2[\overline\chi,\overline\chi],\check\xi\wedge\check\eta>&=&[i(\check\xi)\overline\chi,i(\check\eta)\overline\chi]-i\left([i(\check\xi)\overline\chi,\check\eta]-[i(\check\eta)\overline\chi,\check\xi]-i([\check\xi,\check\eta])\overline\chi \right)\overline\chi\\
&=&[v-\xi_H,w-\eta_H]
-\left(
[v-\xi_H,w]-t_*(i(v-\xi_H)D\eta)\right)\\
&&+\left([w-\eta_H,v]-t_*(i(w-\eta_H)D\xi)
\right)\\
&&+\left( [v,w]-t_*(i(v)D\eta-i(w)D\xi-\firstbra{\xi}{\eta}_\infty)
\right) \\
&=&[\xi_H,\eta_H]-t_*(i(\xi_H)D\eta-i(\eta_H)D\xi+\firstbra{\xi}{\eta}_\infty)\\
&=&[\xi_H,\eta_H]-t_*([\xi,\eta]_\infty)=0.
\end{array}
$$

(iii) It follows from  (\ref{buv}) that, for $\textbf u=\tilde\omega\otimes \tilde{\xi}$, 
$$
\begin{array}{rcl}
[\overline\chi,\textbf u]&=&\mathfrak{L}(\overline\chi)\tilde\omega\otimes\tilde\xi+(-1)^r\tilde\omega\wedge [\overline\chi,\tilde\xi]-(-1)^{2r}d\tilde\omega\wedge i(\tilde\xi)\overline\chi\\
&=&d\tilde\omega\otimes\tilde\xi+(-1)^r\tilde\omega\wedge[\overline\chi,\tilde\xi].
\end{array}
$$
As $D$ is characterized by proposition \ref{propriedadesD}, it is enough to prove $[\overline\chi,\tilde\xi]=\overline D\tilde\xi$. It follows from propositions \ref{colchetek} (ii) and \ref{p31} (ii) that
$$\begin{array}{rcl}
i(\check\eta)[\overline\chi,\tilde\xi]&=&\firstbra{i(\check\eta)\overline\chi}{\tilde\xi}_\infty-i(\firstbra{\check\eta}{\tilde\xi}_\infty)\overline\chi\\
&=&\firstbra{w-\eta_H}{\tilde\xi}_\infty-i(\firstbra{w-\eta_H+\tilde\eta}{\tilde\xi}_\infty)\overline\chi\\
&=&\firstbra{w-\eta_H}{\tilde\xi}_\infty-i(\firstbra{w-\eta_H}{\tilde\xi}_\infty)\overline\chi\\
&=&[w-\eta_H,\xi_H]+i(w-\eta_H)D\xi-([w-\eta_H,\xi_H]-t_*(i(w-\eta_H)D\xi)\\
&=&\rho_*(i(w-\eta_H)D\xi)\\
&=&i(\check\eta)\nu^{-1}D\nu\tilde\xi=i(\check\eta)\overline D\tilde\xi.
\end{array}$$
\EPf

If $\textbf u,\textbf v\in \wedge\tilde{\mathcal T}^*\otimes \tilde{\mathcal J}^{\infty}\mathfrak g$, with $\deg \textbf u=r$, $\deg \textbf v=s$, then we get from (\ref{brauv}) and proposition \ref{p42} (iii) that
\begin{equation}\label{Dcolch}
\overline D[\textbf u,\textbf v]=[\overline D\textbf u,\textbf v]+(-1)^r[\textbf u,\overline D\textbf v],
\end{equation}
and
$$[\overline\chi,[\overline\chi,\textbf u]]=[[\overline\chi,\overline\chi],\textbf u]-[\overline\chi,[\overline\chi,\textbf u]]=-[\overline\chi,[\overline\chi,\textbf u]],$$
therefore, $\overline D^2\textbf u=0$, or
$$\overline D^2=0.$$
Then it is well defined \emph{the second linear Spencer complex},
$$0\rightarrow \Gamma(\mathfrak g)\stackrel{j^\infty}{\rightarrow}\tilde{\mathcal J}^{\infty}\mathfrak g\stackrel{\overline D}{\rightarrow}{\tilde{\mathcal T}}^*\otimes \tilde{\mathcal J}^{\infty}{\mathfrak g}\stackrel{\overline D}{\rightarrow}
{\wedge}^2\tilde{\mathcal T}^*\otimes \tilde{\mathcal J}^{\infty}{\mathfrak g}\stackrel{\overline D}
{\rightarrow}\cdots
\stackrel{\overline D}{\rightarrow}{\wedge}^m\tilde{\mathcal T}^*\otimes\tilde {\mathcal J}^{\infty}\mathfrak g\rightarrow 0,$$
where $\dim T=m$. This complex projects on
$$0\rightarrow \Gamma(\mathfrak g)\stackrel{j^k}{\rightarrow} \tilde{\mathcal J}^{k}\mathfrak g\stackrel{\overline D}{\rightarrow}\tilde{\mathcal T}^*\otimes \tilde{\mathcal J}^{k-1}\mathfrak g\stackrel{\overline D}{\rightarrow}{\wedge^2}\tilde{\mathcal T}^*\otimes \tilde{\mathcal J}^{k-2}\mathfrak g\stackrel{\overline D}{\rightarrow}\cdots 
 \stackrel{\overline D}{\rightarrow}\wedge^{m}\tilde{\mathcal T}^*\otimes \tilde{\mathcal J}^{k-m}\mathfrak g\rightarrow 0.$$
This complex is exact (see \cite {Ma1}, \cite{Ma2}, \cite{KS}).

%\subsection {The second nonlinear Spencer operator}

Let's now introduce the \emph{second nonlinear Spencer operator} $\overline{\mathcal D}$. The ``finite'' form $\overline{\mathcal D}$ of the linear Spencer  operator $\overline D$ is defined by
\begin{equation*}\label{dxi2}
\overline{\mathcal D}\sigma=\overline\chi-\sigma_*^{-1}(\overline\chi),
\end{equation*}
where $\sigma\in\mathcal G_{\infty}$.
\begin{po}\label{overd}
The operator $\overline{\mathcal D}$ take values in $\tilde{\mathcal T}^*\otimes\tilde {\mathcal J}^\infty\mathfrak g$, so
$$\overline{\mathcal D}:\mathcal G_{\infty}\rightarrow \tilde{\mathcal T}^*\otimes \tilde{\mathcal J}^\infty\mathfrak g,$$
and
\begin{equation}\label{ivdsii}
i(v)(\overline{\mathcal D}\sigma)_k=v-
f_*^{-1}(\sigma_{k+1}.v.\sigma_{k+1}^{-1})+(v-j^1\sigma_{k}^{-1}.\sigma_{k+1}.v),
\end{equation}
where $\sigma  = \lim \emph{proj }\sigma_k\in\mathcal G_{\infty}$.
\end{po}
\Pf  Applying (\ref{sigmaomega}) and (\ref{sigmastar}), it follows for  $\tilde{\xi}\in \tilde {\mathcal J}^\infty\mathfrak g$,
$$i(\tilde{\xi})\overline{\mathcal D}\sigma=i(\tilde\xi)\overline{\chi}-\sigma^{-1}_*(i(\sigma_*(\tilde\xi))\overline\chi)=0,$$
and for $v\in{\mathcal T}$,
$$
i(v)\overline{\mathcal D}\sigma=i(v)\overline{\chi}-\sigma^{-1}_*(i(\sigma_*(v))\overline{\chi}),
$$
therefore  from proposition \ref{f41}
$$
\begin{array}{rcl}
i(v)(\overline{\mathcal D}\sigma)_k&=&i(v)\overline{\chi}-(\sigma^{-1}_{k+1})_*\left(i(  f_*(v)+j^1\sigma_k.v.\lambda^1\sigma_{k+1}^{-1}-\lambda^1\sigma_{k+1}.v.\lambda^1\sigma_{k+1}^{-1})      \overline{\chi}\right)
\\&=&v-(\sigma_{k+1}^{-1})_*\left(f_*v-t_*(j^1\sigma_k.v.\sigma_{k+1}^{-1}-\sigma_{k+1}.v.\sigma_{k+1}^{-1})\right)\\
&=&v-(\sigma_{k+1}^{-1})_*\left(f_*v-(j^1\sigma_k.v.j^1\sigma_{k}^{-1}-\sigma_{k+1}.v.\sigma_{k+1}^{-1})\right)\\
&=&v-(\sigma_{k+1}^{-1})_*\left(f_*v-f_*v+\sigma_{k+1}.v.\sigma_{k+1}^{-1}\right)\\
&=&v-(\sigma_{k+1}^{-1})_*\left(\sigma_{k+1}.v.\sigma_{k+1}^{-1}\right)\\
&=&v-\left(
j^1\sigma_k^{-1}.(\sigma_{k+1}.v.\sigma_{k+1}^{-1}).j^1\sigma_k+j^1\sigma_{k}^{-1}.(\sigma_{k+1}.v.\sigma_{k+1}^{-1}).\sigma_{k+1}\right.
\\&& \left. -\sigma_{k+1}^{-1}.(\sigma_{k+1}.v.\sigma_{k+1}^{-1}).\sigma_{k+1}\right)\\
&=&v-
j^1\sigma_k^{-1}.(\sigma_{k+1}.v.\sigma_{k+1}^{-1}).j^1\sigma_k+(v-j^1\sigma_{k}^{-1}.\sigma_{k+1}.v)\\
&=&v-f_*^{-1}(\sigma_{k+1}.v.\sigma_{k+1}^{-1})+(v-j^1\sigma_{k}^{-1}.\sigma_{k+1}.v).\\
\end{array}
$$
Since that
$$
t_*(v-j^1\sigma_{k}^{-1}.\sigma_{k+1}.v)=t_*(v-j^1\sigma_k^{-1}.(\sigma_{k+1}.v.\sigma_{k+1}^{-1}).j^1\sigma_k)=v-
f_*^{-1}(\sigma_{k+1}.v.\sigma_{k+1}^{-1})
$$ we get
%the left action by $j^1\sigma_{k}^{-1}.\sigma_{k+1}$ depends only on $\sigma_{k}^{-1}.\sigma_{k}=I_k$. Then it follows
$$\nu\left(i(v)(\overline{\mathcal D}\sigma)_k\right)=v-j^1\sigma_{k}^{-1}.\sigma_{k+1}.v.$$
\EPf
\begin{co}
We have $\overline{\mathcal D}\sigma=0$ if and only if $\sigma=j^{\infty}(\pi_0\circ\sigma)$, where $\pi_0:\mathcal G_{\infty} \rightarrow \mathcal  G$.
\end{co}
\begin{co}\label{sigmaDii}
If $\sigma_{k+1}\in\mathcal{G}_{k+1}$, then
$$(\sigma_{k+1})_*(v)=\sigma_{k+1}.v.\sigma_{k+1}^{-1}+(\sigma_{k+1})_*(i(v)\overline{\mathcal{D}}\sigma_{k+1}),$$
for $v\in \mathcal{T}$.
\end{co}
\Pf It follows from proposition \ref{actionsigma} (iii) and proposition \ref{overd} that  
$$
\begin{array}{rcl}
(\sigma_{k+1})_*(i(v)\overline{\mathcal{D}}\sigma_{k+1})&=&f_*\left(v-f_*^{-1}(\sigma_{k+1}.v.\sigma_{k+1}^{-1}) \right)+j^1\sigma_{k}.\left(v-j^1\sigma_{k}^{-1}.\sigma_{k+1}.v\right).\sigma_{k}^{-1}\vspace{.2cm}\\
&=&f_*v-\sigma_{k+1}.v.\sigma_{k+1}^{-1} +j^1\sigma_{k}.v.\sigma_{k+1}^{-1}-\sigma_{k+1}.v.\sigma_{k+1}^{-1}\vspace{.2cm}\\
&=&(\sigma_{k+1})_*(v)-\sigma_{k+1}.v.\sigma_{k+1}^{-1}.
\end{array}
$$
\EPf
\begin{po}\label{compo}
The operator $\overline{\mathcal D}$ has the following properties:
\begin{enumerate}
\item [(i)] If $\sigma,\sigma'\in \mathcal G_{\infty}$,
$$\overline{\mathcal D}(\sigma'\circ\sigma)=\overline{\mathcal D}\sigma+\sigma_*^{-1}(\overline{\mathcal D}\sigma').$$
In particular
$$\overline{\mathcal D} \sigma^{-1}=-\sigma_*(\overline{\mathcal D}\sigma).$$
\item [(ii)] If $\sigma\in\mathcal G_{\infty}$, $\mbox{\bf u}\in\wedge\tilde{\mathcal T}^*\otimes \tilde{\mathcal J}^\infty\mathfrak g$,
$$ \overline{\mathcal D}(\sigma^{-1}_*\mbox{\bf u})=\sigma_*^{-1}(\overline{D}\mbox{\bf u})+[\overline{\mathcal D}\sigma,\sigma^{-1}_*\mbox{\bf u}].$$
\item [(iii)] If $\xi=\frac{d}{du}\sigma_u|_{u=0}$, with $\xi\in {\cal J}^\infty\mathfrak g$, and $\sigma_t\in\mathcal G_{\infty}$ is the 1-parameter group associated to $\xi$, then
$$\overline{D}\tilde\xi=\frac{d}{du}\overline{\mathcal D}\sigma_u|_{u=0}.$$
\end{enumerate}
\end{po}
\Pf 
(i) $$\overline{\mathcal D}(\sigma'\circ\sigma)=\overline{\chi}-\sigma_*^{-1}(\overline{\chi})+\sigma_*^{-1}(\overline{\chi}-(\sigma')_*^{-1}(\overline{\chi}))=\overline{\mathcal D}\sigma+\sigma_*^{-1}(\overline{\mathcal D}\sigma').$$

(ii) $$\begin{array}{rcl}
 \overline{D}(\sigma_*^{-1}\textbf u)&=&\sigma_*^{-1}[\sigma_*\overline{\chi},\textbf u]=\sigma_*^{-1}[\overline{\chi}-\overline{\mathcal D}\sigma^{-1},\textbf u]\vspace{.3cm}\\
&=&\sigma_*^{-1}(\overline{D}\textbf u)-[\sigma_*^{-1}(\overline{\mathcal D}\sigma^{-1}),\sigma_*^{-1}\textbf u]=\sigma_*^{-1}(\overline{D}\textbf u)+[\overline{\mathcal D}\sigma,\sigma_*^{-1}\textbf u].
\end{array}$$

(iii) $$\frac{d}{du}\overline{\mathcal D}\sigma_u|_{u=0}=-\frac{d}{du}(\sigma_u^{-1})_*(\overline{\chi})=-[\xi,\overline{\chi}]=-[\tilde\xi-\xi_H,\overline{\chi}]=\overline{D}\tilde\xi+[\xi_H,\overline{\chi}].$$
It follows from Proposition \ref{p31}(ii) that
$$
\begin{array}{rcl}
i(\check\eta)[\overline{\chi},\xi_H]&=&\firstbra{i(\check\eta)\overline{\chi}}{\xi_H}_\infty-i(\firstbra{\check\eta}{\xi_H}_\infty)\overline{\chi}\\
&=&\firstbra{v-\eta_H}{\xi_H}_\infty-i(\firstbra{\tilde\eta+v-\eta_H}{\xi_H}_\infty)\overline{\chi}\\
&=&\firstbra{v-\eta_H}{\xi_H}_\infty-i(\firstbra{v-\eta_H}{\xi_H}_\infty)\overline{\chi}-i(\firstbra{\tilde\eta}{\xi_H}_\infty)\overline{\chi}\\
&=&i(\firstbra{\xi_H}{\tilde\eta}_\infty)\overline{\chi},
\end{array}
$$
where $\check\eta=v+\eta$, $\eta_H=t_*\eta$, $\xi_H=t_*\xi$.
See that $i(\firstbra{\tilde\eta}{\xi_H}_\infty)\overline{\chi}=0$.
\EPf

Proposition \ref{overd} says that $\overline{\mathcal D}$ is projectable:
\begin{equation*}\label{f42ii}
\begin{array}{rcl}
\overline{\mathcal D}:\mathcal G_{k+1}&\rightarrow &\tilde{\mathcal T}^*\otimes \tilde{\mathcal J}^k\mathfrak g\\
\sigma_{k+1}&\mapsto &\overline{\mathcal D}\sigma_{k+1}
\end{array},
\end{equation*}
where 
\begin{equation}\label{bard}
i(v)\overline{\mathcal D}\sigma_{k+1}=v-
f_*^{-1}(\sigma_{k+1}.v.\sigma_{k+1}^{-1})+(v-j^1\sigma_{k}^{-1}.\sigma_{k+1}.v).
\end{equation}
%(see formula (6.8) of \cite{Ma2})

It follows from $[\overline{\chi},\overline{\chi}]=0$ that 
$$0=\sigma_*^{-1}([\overline{\chi},\overline{\chi}])=[\sigma_*^{-1}(\overline{\chi}),\sigma_*^{-1}(\overline{\chi})]=[\overline{\chi}-\overline{\mathcal D}\sigma,\overline{\chi}-\overline{\mathcal D}\sigma]=[\overline{\mathcal D}\sigma,\overline{\mathcal D}\sigma]-2\overline{D}(\overline{\mathcal D}\sigma),$$
therefore
\begin{equation}\label {fs43}
\overline{D}(\overline{\mathcal D}\sigma)-\frac{1}{2}[\overline{\mathcal D}\sigma,\overline{\mathcal D}\sigma]=0.
\end{equation}
If we define the non linear operator 
\begin{equation*}\label{fs44}
\begin{array}{rcl}
\overline{\mathcal D}_1:\tilde{\mathcal T}^*\otimes \tilde{J}^\infty\mathfrak g&\rightarrow &\wedge^2\tilde{\mathcal T}^*\otimes \tilde{J}^\infty{\mathfrak g}\\
\textbf u&\mapsto & \overline{D}\textbf u-\frac{1}{2}[\textbf u,\textbf u]
\end{array},
\end{equation*}
then we can write  (\ref{fs43}) as
$$\overline{\mathcal D}_1\overline{\mathcal D}=0.$$
The operator $\overline{\mathcal D}_1$ projects in order $k$ to
$$\overline{\mathcal D}_1:\tilde{\mathcal T}^*\otimes \tilde{\mathcal J}^k\mathfrak g\rightarrow \wedge^2\tilde{\mathcal T}^*\otimes \tilde{\mathcal J}^{k-1}\mathfrak g,$$
where
$$\overline{\mathcal D}_1\textbf u=\overline{D}\textbf u-\frac{1}{2}[\textbf u,\textbf u]_k.$$
We define the \emph{second non-linear Spencer complex} by
$$1\rightarrow {\mathcal G}\stackrel{j^{k+1}}{\rightarrow}\mathcal G_{k+1}\stackrel{\overline{\mathcal D}}{\rightarrow} \tilde{\mathcal T}^*\otimes \tilde{\mathcal J}^k\mathfrak g\stackrel{\overline{\mathcal D}_1}{\rightarrow}\wedge^2\tilde{\mathcal T}^*\otimes \tilde{\mathcal J}^{k-1}\mathfrak g,$$
which is exact in $\mathcal G_{k+1}$.

Let be $p^1(\overline{\mathcal D}):J^1G_{k+1}\rightarrow \tilde T^*\otimes\tilde J^k\mathfrak g$ the morphism associated to the differential operator $\overline{\mathcal D}$. It follows from (\ref{ivdsii}) that 
\begin{equation*}%\label{p1d}
p^1(\overline{\mathcal D})(j^1_x\sigma_{k+1})(v)=t_*(v-j^1_x\sigma_k^{-1}.\sigma_{k+1}(x).v)+v-j^1_x\sigma_k^{-1}.\sigma_{k+1}(x).v
\end{equation*}
where $\sigma_{k+1}$ is an admissible section of $G_{k+1}$ and $\sigma_k=\pi_k\sigma_{k+1}$.

\begin{po} The image of $G_{k,1}$ by $p^1(\overline{\mathcal D})$ is the set
$$
\tilde B^{k,1}=\{X\in {\tilde T}^*\otimes \tilde J^k\mathfrak g: v\in T\rightarrow t_*(v-\nu(X(v)))\in T \mbox{ is inversible}\}.
$$
\end{po}
\Pf It follows from (\ref{bard}) that 
$$
t_*(\nu(i(v)\overline{\mathcal D}\sigma_{k+1}))=v-
f_*^{-1}(\sigma_{k+1}.v.\sigma_{k+1}^{-1}),
$$
 so $t_*(v-\nu(i(v)\overline{\mathcal D}\sigma_{k+1}))=f_*^{-1}(\sigma_{k+1}.v.\sigma_{k+1}^{-1})$ therefore $\overline{\mathcal D}\sigma_{k+1}\in \tilde B^{k,1}$. 

Conversely, let be $Y\in \tilde B^{k,1}_x$ and consider $\bar Y\in \tilde B^{k+1,1}_x$ such that $\mbox{id}^*\otimes\pi_k(\bar Y)=Y$. Let be the map $\Sigma:T_x\rightarrow T_{I_{k+1}(x)}G_k$ defined by $\Sigma(v)=v-\nu(i(v)\bar Y)$.   There exists a local section 
$\sigma\in \mathcal{G}_{k+1}$ such that $\sigma(x)=I_{k+1}(x)$ and $\sigma_*^{-1}(v)=\Sigma(v)$. Then 
$$
\begin{array}{rcl}
i(v)p^1(\overline{\mathcal D})(j^1_x\sigma)&=&t_*(v-j^1_x(\pi_k\sigma)^{-1}.\sigma(x).v)+v-j^1_x(\pi_k\sigma)^{-1}.\sigma(x).v\\
&=&t_*(v-j^1_x(\pi_k\sigma)^{-1}.v)+v-j^1_x(\pi_k\sigma)^{-1}.v\\
&=&t_*(v-\pi_ki(v)\Sigma)+v-\pi_ki(v)\Sigma\\
&=&t_*(v-(v-\nu(i(v) Y))+v-(v-\nu(i(v) Y)\\
&=&i(v)Y.
\end{array}
$$
\EPf

\section{The sophisticated Spencer complex}\label{sSc}

Let be $\overline\delta$ the restriction of $-\overline D$ to $\wedge\tilde T\otimes \gamma^k$. Set
$$
B^{k,p}=\wedge^p\tilde T^*\otimes \tilde J^k\mathfrak g\,\,/\,\,
%\setminus
\overline\delta\left(\wedge^{p-1} \tilde T^*\otimes  \gamma^{k+1}\right)
$$
and
$$
B^k=\oplus B^{k,p}. 
$$
Let $\mathcal{B}^{k}$ be the sheaf of sections of $B^k$. 
\begin{po}
The subsheaf $\Gamma\left(\overline\delta(\wedge\tilde T^*\otimes \mathfrak \gamma^{k+1})\right)$ is an ideal of $\wedge\tilde{\mathcal{T}}^*\otimes\tilde{\mathcal{J}}^k\mathfrak g$.
\end{po}
\Pf 
Let be $\textbf u\in \wedge\tilde{\mathcal T}^*\otimes \tilde{\mathcal J}^{k}\mathfrak g$, with $\deg \textbf u=r$, $\textbf v\in\Gamma(\wedge \tilde T^*\otimes \gamma^{k+1})$  with $\deg \textbf v=s-1$ and   $\textbf u_{k+1}\in \wedge\tilde{\mathcal T}^*\otimes \tilde{\mathcal J}^{k+1}\mathfrak g$ such that $\pi_k\textbf u_{k+1}=\textbf u_{k}$. 
Then we get from $\pi_k[\textbf u_{k+1},\textbf v]_{k+1}=[\textbf u,\pi_{k}\textbf v]_{k}=0$ that $[\textbf u_{k+1},\textbf v]_{k+1}\in \Gamma(\wedge \tilde T^*\otimes \gamma^{k+1})$.
Therefore we get from (\ref{Dcolch}) that
$$-\overline \delta[\textbf u_{k+1},\textbf v]_{k+1}=[\overline D\textbf u_{k+1},\pi_{k}\textbf v]_k+(-1)^r[\textbf u,-\overline \delta\textbf v]_k=(-1)^r[\textbf u,-\overline \delta\textbf v]_k$$
and the proposition is proved.
\EPf

\begin{co}
$\mathcal{B}^{k}$ is a Lie algebra sheaf. 
\end{co}
%$\mbox{ incorporar o ethos do morador urbano economicamente desfavorecido }$

Given $\mathbf u\in {\mathcal B}^{k,p}$ let be $\mathbf u_{k+1}\in \wedge^p\tilde{\mathcal T}^*\otimes\tilde{\mathcal J}^{k+1}\mathfrak g$ that projects on $\mathbf u$. We define the operator 
$$
\hat D:{\mathcal B}^{k,p}\rightarrow {\mathcal B}^{k,p+1}
$$ 
by 
$$
\hat D\mathbf u=\overline D\mathbf u_{k+1}\mod \overline \delta (\wedge^{p}{\tilde{\mathcal T}}^*\otimes \gamma^{k+1})
$$
The sophisticated Spencer complex is
\begin{equation}\label{linsoph}
0\rightarrow j^k\mathcal G\rightarrow {\mathcal B}^{k}\stackrel{\hat{D}}{\rightarrow} {\mathcal B}^{k,1}\stackrel{\hat{D}}{\rightarrow}{\mathcal B}^{k,2}\stackrel{\hat{D}}{\rightarrow}\cdots\stackrel{\hat{D}}{\rightarrow} {\mathcal B}^{k,p}\stackrel{\hat{D}}{\rightarrow} {\mathcal B}^{k,p+1}\stackrel{\hat{D}}{\rightarrow}\cdots\stackrel{\hat{D}}{\rightarrow} {\mathcal B}^{k,n}\rightarrow 0
\end{equation}
which is exact

We introduce now the nonlinear version of (\ref{linsoph}).  Let be 
$$
G^{k+1}_k(x)=\{X\in G_{k+1}(x)|\pi^{k+1}_kX=I_k(x)\}.
$$
We know that $J^1G_k$ has a affine structure on $G_k$. The vector space for this affine structure is the fiber of
$$
vJ^1G_k=\{w\in TJ^1G_k:\pi^1_0w=0\}.
$$
Therefore we can identify 
$$\begin{array}{ccl}
\partial: G^{k+1}_k (x)&\rightarrow & T^*_x\otimes \gamma_x^{k}\\
X&\rightarrow & \partial X=X-I_{k+1}(x)
\end{array}
$$
It follows from (\ref{bard}) that for a section $F\in \mathcal G^{k+1}_k$ we get $\overline{\mathcal D} F=-\overline \delta (\nu_k^{-1}\partial F)$, where $\partial F=F-(I_{k+1})$.

Let be $F\in{\mathcal G}_k$ and $F_1\in \mathcal G_{k+1}$ such that $\pi_kF_1=F$ and let be 
$$
\widehat{\mathcal D}F=\bar{\mathcal D}F_1\mod\bar\delta(\gamma^{k+1})\in  {\mathcal B}^{k,1}
$$
Let's prove that this class depends only of $F$. If $F_2$ is another section in $ \mathcal G_{k+1}$ such that $\pi_kF_2=F$, then $F_2=F_1G$ with $G\in\mathcal G_{k}^{k+1}$. Then
$$
\overline{\mathcal D}F_2=\overline{\mathcal D}(F_1G)=\overline{\mathcal D}G+G_*^{-1}(\overline{\mathcal D}F_1).
$$
It follows from proposition \ref{actionsigma} (iii) that $G^{-1}$ acts trivially on $\overline{\mathcal D}F_1$. Therefore
$$
\overline{\mathcal D}F_2=-\overline\delta(\nu_k^{-1}\partial G)+\overline{\mathcal D}F_1,
$$
what proves the class depends only of $F$.

Let's define for $\textbf u\in {\mathcal B}^{k,1}$ 
$$
\hat{\mathcal D}_1\textbf u=\hat D\textbf u-\frac 1 2[\textbf u,\textbf u]
$$
We obtain the nonlinear sophisticated complex
$$
1\rightarrow {\mathcal G}\stackrel{j^{k+1}}{\rightarrow}\mathcal G_{k+1}\stackrel{\hat{\mathcal D}}{\rightarrow} {\mathcal B}^{k,1}\stackrel{\hat{\mathcal D}_1}{\rightarrow}{\mathcal B}^{k,2},
$$
which is exact in $\mathcal G_{k+1}$.

%XXXXX
%Let be $p^1(\overline{\mathcal D}):J^1G_{k+1}\rightarrow \tilde T^*\otimes\tilde J^k\mathfrak g$ the morphism associated to the differential operator $\overline{\mathcal D}$. It follows from proposition \ref{ivdsii} that 
%\begin{equation}\label{p1d}
%p^1(\overline{\mathcal D})(j^1_x\sigma_{k+1})(v)=t_*(v-j^1_x\sigma_k^{-1}.\sigma_{k+1}(x).v)+v-j^1_x\sigma_k^{-1}.\sigma_{k+1}(x).v
%\end{equation}
%where $\sigma_{k+1}$ is an admissible section of $G_{k+1}$ and $\sigma_k=\pi_k\sigma_{k+1}$.
%XXXXXXX

%If $p^1(\hat{\mathcal{D}}):J^1G_{k+1}\rightarrow B^{k,1}$ is the morphism associated to the operator   $\hat{\mathcal{D}}$, it follows from (\ref{p1d}) that
%\begin{equation*}
%p^1(\hat{\mathcal D})(j^1_x\sigma_{k+1})(v)=t_*(v-j^1_x\sigma_k^{-1}.\sigma_{k+1}(x).v)+v-j^1_x\sigma_k^{-1}.\sigma_{k+1}(x).v\mod\bar\delta(S^{k+1}\mathcal T^*\otimes\Gamma(\mathfrak g))
%\end{equation*}

%Consider $\mbox{id}^*\otimes t_* :\tilde T^*\otimes \tilde J^k\mathfrak g\rightarrow \tilde T^*\otimes \tilde T$. 
%Let be $\widehat B^{k,1}$ the set of ${\textbf u}\in B^{k,1}$ such that $\mbox{id}^*\otimes\mbox{id}-(\mbox{id}^*\otimes t_*)\textbf{u}$ is inversible.

\end{document}